\documentclass[10pt,twoside]{siamltex}
\usepackage{amsfonts,epsfig}

\usepackage{graphicx}
\usepackage{amssymb}
\usepackage{amsmath}
\usepackage{enumerate}

\usepackage{color}

\usepackage{eqparbox}

\usepackage{algorithm}
\usepackage{algorithmic}

\usepackage{xcolor}
\usepackage{tikz}

\definecolor{mygray}{RGB}{47,79,79}
\newcommand{\tc}{\textcolor{black}}

\newcommand{\N}{\mathbb N}
\newcommand{\R}{\mathbb R}

\newcommand{\V}{\mathbb V}

\newcommand{\dist}{\operatorname{dist}}
\newcommand{\Exp}{\operatorname{Exp}}
\newcommand{\Log}{\operatorname{Log}}

\newcommand{\id}{\operatorname{id}}

\makeatletter
\renewcommand*\env@matrix[1][*\c@MaxMatrixCols c]{%
  \hskip -\arraycolsep
  \let\@ifnextchar\new@ifnextchar
  \array{#1}}
\makeatother

\newcommand{\mcM}{\mathcal{M}}
\newcommand{\mcN}{\mathcal{N}}

\newtheorem{remark}{Remark}

\title{Hermite Interpolation and data processing errors on Riemannian Matrix Manifolds}
\author{
  Ralf Zimmermann\thanks{Department of Mathematics and Computer Science, University of Southern Denmark (SDU) Odense,
    (zimmermann@imada.sdu.dk).}
}
\begin{document}

\maketitle

\begin{abstract}
The main contribution of this paper is twofold:
On the one hand, a general framework
for performing Hermite interpolation on Riemannian manifolds is presented. The method is applicable, if algorithms
for the associated Riemannian exponential and logarithm mappings are available. This includes many of the matrix manifolds
that arise in practical Riemannian computing application such as data analysis and signal processing, computer vision
and image processing, structured matrix optimization problems
and model reduction.

On the other hand, we expose a natural relation between data processing errors and the sectional curvature of the manifold in question.
This provides general error bounds for manifold data processing methods that rely on Riemannian normal coordinates.

Numerical experiments are conducted for the compact Stiefel manifold
of rectangular column-orthogonal matrices.
As use cases, we compute Hermite interpolation curves for orthogonal
matrix factorizations such as the 
singular value decomposition and the QR-decomposition.
%
%
\end{abstract}



\section{Introduction}
Given a data set that consists of locations $t_0,\ldots, t_k\in \R$, function values $f_0 = f(t_0),\ldots, f_k = f(t_k)$
and derivatives $\dot{f}_0 = \dot{f}(t_0),\ldots, \dot{f}_k=\dot{f}(t_k)$,
the (first-order) Hermite interpolation problem reads:
\begin{quote}
Find a polynomial $P$ of suitable degree such that 
 \begin{equation}
 \label{eq:basicHermite}
 P(t_i) = f_i, \quad \dot{P}(t_i) = \dot{f}_i,  \quad i =0,\ldots,k.
\end{equation}
\end{quote}
{\em Local cubic Hermite interpolation} is the special case of Hermite-interpolating a two-points data set
$\{f_{i}, \dot{f}_{i}, f_{i+1}, \dot{f}_{i+1}\}$ on $t_i,t_{i+1}\in \R$.
{\em Cubic Hermite interpolation} is achieved by joining the local pieces on each sub-interval $[t_i,t_{i+1}]$.
By construction, the derivative at the end point of $[t_i,t_{i+1}]$ coincides with the derivative
of the start point of $[t_{i+1},t_{i+2}]$ so that the resulting curve is globally $C^1$, 
\cite[Remark 7.7]{hohmann2003numerical}.

In this paper, we address the Hermite interpolation problem for a function that takes values on a Riemannian manifold $\mcM$
with tangent bundle $T\mcM$. More precisely, consider a differentiable function
\[
 f:[a,b] \to \mcM, \quad t\mapsto f(t)
\]
and a sample plan $a=t_0,\ldots, t_k=b$.
Sampling of the function values and the derivatives of $f$ at the parameter instants $t_i$
produces a data set consisting of manifold locations $p_i = f(t_i)\in \mcM$
and velocity vectors $v_{p_i}\in T_{p_i}\mcM$ in the respective tangent spaces of $\mcM$ at $p_i$.
The Hermite manifold interpolation problem is:
\begin{quote}
Find a curve $c:[a,b] \to \mcM$ of class $C^1$ such that
 \begin{equation}
 \label{eq:basicHermiteMnf}
 c(t_i) = p_i\in \mcM, \quad \dot{c}(t_i) = v_{p_i}\in T_{p_i}\mcM,  \quad i =0,\ldots,k.
\end{equation}
\end{quote}
\newpage
\subsection{Original contributions}
(1) We introduce a method to tackle problem \eqref{eq:basicHermiteMnf} that is a direct analogue to Hermite interpolation in Euclidean spaces.
The method has the following features:
\begin{enumerate}[(i)]
 \item The approach works on arbitrary Riemannian manifolds, i.e., no special structure
 (Lie Group, homogeneous space, symmetric space,...) is required.\\
 In order to conduct practical computations, only algorithms for evaluating the Riemannian exponential map
 and the Riemannian logarithm map must be available.\footnote{The Riemannian exp and log maps for some of the 
 most prominent matrix manifolds are collected in \cite{Zimmermann_MORHB2019}.}
 \item The computational effort, in particular, the number of
 Riemannian exp and log evaluations is lower than that of any other Hermite manifold interpolation method known to the author.
\end{enumerate}
(2) In addition, we expose a natural relation between data processing errors and the sectional curvature of the manifold in question.
This provides general error bounds for data processing methods (including but not limited to interpolation) that work via a back-and-forth mapping of data between the manifold and its tangent space, or, more precisely, data processing methods that rely on Riemannian normal coordinates.

For convenience, the exposition will focus on cubic polynomial Hermite interpolation. However, the techniques may be readily combined with any interpolation method that is linear in the sampled locations and derivative values.
Apart from polynomial interpolation, this includes radial basis function approaches \cite{Amsallem2010} and gradient-enhanced Kriging \cite{Zimmermann2013_GEK}.

As a use-case, we provide an explicit and efficient method for the cubic Hermite interpolation of column-orthogonal matrices, which form the so-called Stiefel manifold $St(n,r) = \{U\in\R^{n\times r}| U^TU=I\}$.
Stiefel matrices arise in orthogonal matrix factorizations such as the singular value decomposition
and the QR-decomposition.
\subsection{Related work}
Interpolation problems with ma\-ni\-fold-valued sample data and spline-related approaches have triggered an extensive amount of research work.

It is well-known that cubic splines in Euclidean spaces are acceleration-minimi\-zing.
This property allows for a generalization to Riemannian manifolds in form of a variational problem
for the intrinsic, covariant acceleration of curves,
whose solutions can be interpreted as generalized cubic polynomials on Riemannian manifolds.
The variational approach to interpolation on manifolds has been investigated e.g. in
\cite{Noakes1989, Crouch1995, CAMARINHA2001107, Steinke2008,BOUMAL2011,Samir2012, Kim2018},
see also \cite{Noakes2007} and references therein.
While the property of minimal mean-acceleration is certainly desirable in many a context, including automobile, aircraft and ship designs and digital animations,
there is no conceptual reason to impose this condition when interpolating general smooth non-linear manifold-valued functions.

\begin{figure}[ht]
\centering
\includegraphics[width=0.4\textwidth]{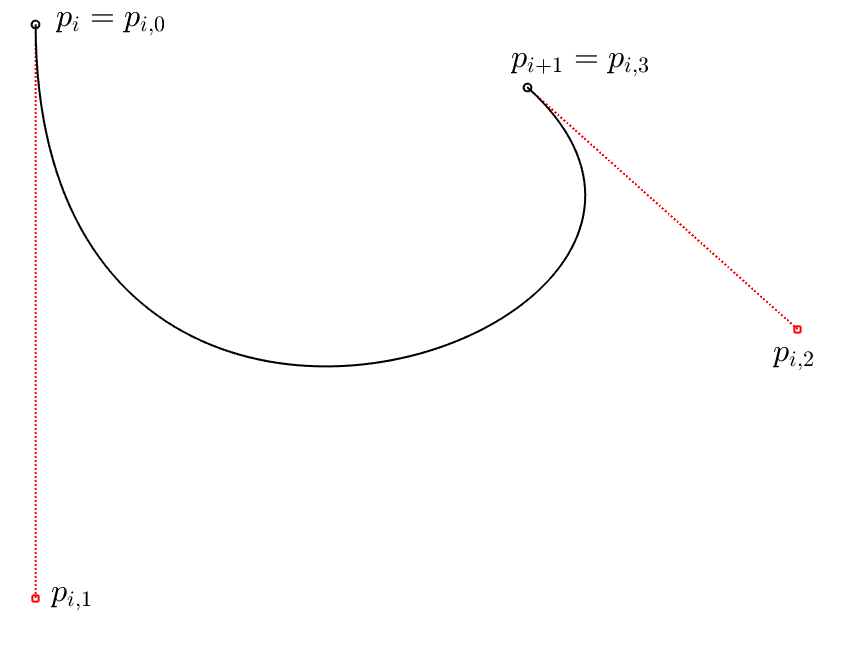}
\caption{A cubic B{\'e}zier curve based on four control points $p_{i,0},p_{i,1},p_{i,2}, p_{i,3}$.
The `inner' control points $p_{i,1},p_{i,2}$ may be used to prescribe tangent directions at $p_i=p_{i,0}$  and 
$p_{i+1}=p_{i,3}$, which are interpolated.
}
\label{fig:CubicBezier}
\end{figure}

%
A related line of research is the generalization of B{\'e}zier curves and the De Casteljau-algorithm \cite{bartels1995introduction}
to Riemannian manifolds \cite{Noakes2007, Krakowski2015SOLVINGIP, Polthier2013, AbsGouseWirth2016, GouseMassartAbsil2018, SAMIR2019}.
B{\'e}zier curves in Euclidean spaces are polynomial splines that rely on a number of so-called control points.
A B{\'e}zier curve starts at the first control point and ends at the last control point,
the starting velocity is tangent to the line between the first two-pair of control points;
the velocity at the endpoint is tangent to the line between the penultimate and the last control point.
This is illustrated in Fig. \ref{fig:CubicBezier}.
The number of control points determines the degree of the polynomial spline.
To obtain the value $B(t)$ of a B{\'e}zier curve at time $t$, a recursive sequence of straight-line convex combinations
of two locations must be computed.
The transition of this technique to Riemannian manifolds is via replacing the inherent straight lines with geodesics \cite{Noakes2007}.
The start and end velocities of the resulting spline are proportional to the velocity vectors of 
the geodesics that connect the first two and the last two control points, respectively \cite[Theorem 1]{Noakes2007}.
%
%

Note that the actual applications and use cases featured in the work referenced above are almost exclusively on low-dimensional matrix manifolds like $S^2, S^3, SO(3)$ or $SE(3)$.

%

A Hermite-type method that is specifically tailored for interpolation problems on the Grassmann manifold is sketched in \cite[\S 3.7.4]{Amsallem2010}.
General Hermitian manifold interpolation has been considered explicitly in \cite{Jakubiak2006}.
The idea is as follows: Given two points $p,q\in\mcM$ on a manifold 
and two tangent directions $v_p\in T_p\mcM, v_q\in T_q\mcM$, the the authors of \cite{Jakubiak2006}
approach the task to construct a connecting curve
$c:[t_i,t_{i+1}]\to \mcM$ such that $c(t_i) = p, \dot{c}(t_i) = v_p, c(t_{i+1}) = q, \dot{c}(t_{i+1}) = v_q$
by constructing a ``left'' arc $l_i$ that starts at $t=t_i$ from $p$ with the prescribed velocity $v_p$
and a ``right'' arc $r_i$ that ends at $t=t_{i+1}$ at $q$ with the prescribed velocity $v_q$.
The two arcs are then blended to a single spline arc via a certain geometric convex combination.
In Euclidean spaces, this would read $s(t) = (1-\Phi(t))l_i(t) + \Phi(t)r_i(t)$,
where $\Phi$ is a suitable weight function.
Because a general Riemannian manifold lacks a vector space structure, the challenge is to 
construct a manifold analogue of a convex combination and \cite{Jakubiak2006} proposes a method
that works on compact, connected Lie groups with a bi-invariant metric.

This same idea of blending a left and a right arc has been followed up in \cite{GouseMassartAbsil2018}.
Here, the Euclidean convex combination is replaced with a geodesics average
$s(t) = \Exp_{l_i(t)}(\Phi(t) \Log_{l_i(t)}(r_i(t)))$.
In combination, this constitutes a valid approach for solving \eqref{eq:basicHermiteMnf}
in arbitrary Riemannian manifolds.\footnote{In practice, the building arcs $l_i(t)$ and $r_i(t)$ may be taken to be the geodesics with the prescribed velocities in their respective start and end points. 
}

It should be mentioned that none of the papers on B{\'e}zier curves referenced above tackle the Hermite interpolation problem explicitly. However, the B{\'e}zier approach can be turned into an Hermite method by choosing the control points such that the sampled start and terminal velocities are met. It is clear that this requires at least $4$ control points in each subinterval $[t_i,t_{i+1}]$, see Fig. \ref{fig:CubicBezier}.

Interpolation problems on Stiefel Manifolds have been considered in \cite{Krakowski2015SOLVINGIP}, however with using quasi-geodesics rather than geodesics.
The work \cite{Zimmermann2018} includes preliminary numerical experiments for interpolating orthogonal frames on the Stiefel manifold that 
relies the canonical Riemannian Stiefel logarithm \cite{Rentmeesters2013, StiefelLog_Zimmermann2017}.

{\bf Remark:} (Hermite) interpolation of curves on Riemannian manifolds, i.e., of manifold-valued functions $f:[a,b] \rightarrow \mcM$
must not be confused with (Hermite) interpolation of real-valued functions with domain of definition on a manifold,
$f:\mcM \rightarrow \R$.
The latter line of research is pursued, e.g., in \cite{Narcowich1995} but is not considered here.
%
%
%
\subsection{Organization}
The paper is organized as follows: Starting from the classical Euclidean case, Section \ref{sec:HermiteInterpMnf}
introduces an elementary approach to Hermite interpolation on general Riemannian manifolds.
Section \ref{sec:errorCurvature} relates the data processing errors of calculations in Riemannian normal coordinates 
to the curvature of the manifold in question. In Section \ref{sec:Stiefel_Hermite}, the specifics of performing Hermite 
interpolation of column-orthogonal matrices are discussed and Section \ref{sec:experiments} illustrates the theory by means 
of numerical examples. Conclusions are given in Section \ref{sec:conclusions}.
\subsection{Notational specifics}
\label{sec:Notation}
The $(r\times r)$-{\em identity matrix} is denoted by $I_r\in\R^{r\times r}$,
or simply $I$ if the dimensions are clear.
The $(r\times r)$-{\em orthogonal group} is denoted by
\[
  O(r) = \{\Phi \in \R^{r\times r}| \Phi^T\Phi = \Phi\Phi^T = I_r\}.
\]
Throughout, the QR-decomposition $A=QR$ of $A\in\R^{n\times r}$, $n\geq r$,
is understood as the `economy size' QR-decomposition with $Q\in \R^{n\times r}$, $R\in \R^{r\times r}$.

The standard matrix exponential and the principal matrix logarithm are defined by
\[
 \exp_m(X):=\sum_{j=0}^\infty{\frac{X^j}{j!}}, \quad \log_m(I+X):=\sum_{j=1}^\infty{(-1)^{j+1}\frac{X^j}{j}}.
\]
The latter is well-defined for matrices that have no eigenvalues on $\R^-$.

\begin{figure}[ht]
\centering
\includegraphics[width=0.7\textwidth]{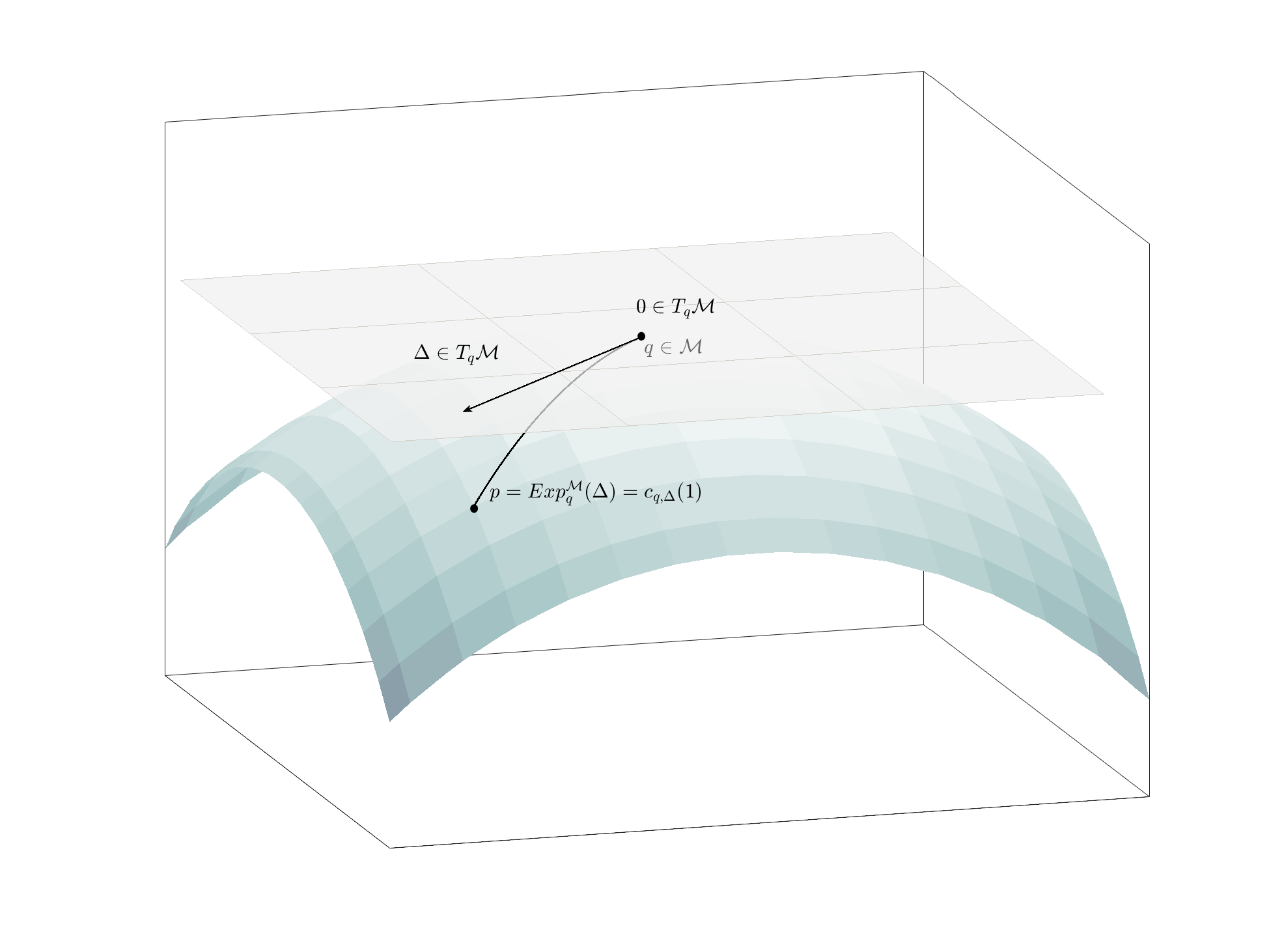}
\caption{Visualization of the Riemannian exponential map: The tangent velocity $\Delta\in T_q\mcM$ is
mapped to the endpoint of the geodesic $p=c_{q,\Delta}(1)\in \mcM$.
The Riemannian distance $\dist(q,p)$ equals the norm $\|\Delta\|$ as measured by the Riemannian metric on $T_q\mcM$.
}
\label{fig:Manifold_plot_RiemannExp}
\end{figure}

For a  Riemannian manifold $\mcM$, the geodesic that starts from $p\in \mcM$ with velocity $v\in T_p\mcM$
is denoted by $t\mapsto c_{p,v}(t)$.
The Riemannian exponential function at $p$ is
\begin{equation}
 \label{eq:RiemannExp}
 \Exp_p^\mcM: T_p\mcM\supset D_0 \rightarrow \mathcal{D}_p\subset \mcM, \quad v\mapsto \Exp_p^\mcM(v):= c_{p,v}(1)
\end{equation}
and maps a small star-shaped domain $D_0$ around the origin in $T_p\mcM$ 
diffeomorphically to a domain $\mathcal{D}_p\subset \mcM$, see Fig. \ref{fig:Manifold_plot_RiemannExp}.
The Riemannian logarithm at $p$ is
\begin{equation}
 \label{eq:RiemannLog}
 \Log_p^\mcM:  \mcM\supset \mathcal{D}_p \rightarrow  D_0\subset T_p\mcM, \quad q\mapsto v = (\Exp_p^\mcM)^{-1}(q).
\end{equation}
Recall that for a differentiable function $f:\mcM\to \mcN$, the differential at $p$ is a 
linear map between the tangent spaces
\begin{equation}
 \label{eq:diff_map}
 df_p: T_p\mcM \to T_{f(p)}\mcN.
\end{equation}
\section{Hermite interpolation on Riemannian manifolds}
\label{sec:HermiteInterpMnf}
In this section, we construct a quasi-cubic spline between two data points $p_0, p_{1}\in \mcM$
on a manifold with prescribed velocities $v_{p_0}\in T_{p_0}\mcM$ and $v_{p_1}\in T_{p_1}\mcM$.
To this end, we develop a manifold equivalent to the classical local cubic Hermite interpolation 
in Euclidean spaces \cite[\S 7]{hohmann2003numerical}.
%
\subsection{The Euclidean case}
\label{sec:HermiteInterpEuclid}
We start with a short recap of Hermite cubic space curve interpolation, where the following setting is of special interest to our considerations.
Let $\V$ be a real vector space and let $f:[t_0,t_1] \to \V$ be differentiable
with $f(t_0) = p\in\V$, $f(t_1) = q=0\in\V$, and derivative data
$v_0 = \dot{f}(t_0), v_1 = \dot{f}(t_1)$.
When applied to vector-valued functions, the classical local cubic Hermite interpolating spline
is the space curve $c(t)$ that is obtained via a linear combination of the sampled data,\footnote{It is an elementary, yet often overlooked fact that for functions $t\mapsto f(t)=(f_1(t),\ldots,f_n(t))^T\in \R^n$,
component-wise polynomial interpolation of the coordinate functions $f_l(t)$ is equivalent 
to interpolating the coefficients in a linear combination of the sampled data vectors.}
\begin{equation}
\label{eq:VectorCubicHermite}
 c(t) = a_0(t)p + a_1(t)q + b_0(t)v_0+ b_1(t)v_1.
\end{equation}
For the reader's convenience, the basic cubic Hermite polynomials coefficient polynomials $a_0(t), a_1(t), b_0(t), b_1(t)$
are listed in Appendix \ref{app:classicHermite}.
%
\subsection{Transfer to the manifold setting}
\label{sec:HermiteInterp2Mnf}
Let $\mcM$ be a Riemannian manifold and consider a differentiable function
\[
 f:[t_0,t_1] \to \mcM, \quad t\mapsto f(t).
\]
Suppose that $f(t_0) = p,f(t_1)=q\in \mcM$ and $\dot f(t_0) = v_p \in T_p\mcM, \dot f(t_1) = v_q \in T_q\mcM$
and assume further that $\dist(p,q)< r_\mcM(q)$, where $r_\mcM(q)$ is the injectivity radius of $\mcM$ at $q$.
The latter condition ensures that the sample data lies within a domain, where the Riemannian normal coordinates are one-to-one, \cite[p. 271]{DoCarmo2013riemannian}.

Our approach is to express the interpolating curve in terms of normal coordinates centered at $q=f(t_1)\in \mcM$,
\[
 c:[t_0,t_1] \rightarrow \mcM, \quad c(t) = \Exp^\mcM_q(\gamma(t)).
\]
%
Hence, the task is transferred to constructing a curve $t\mapsto \gamma(t)\subset T_q\mcM$ 
such that the image curve $c$ under the exponential function solves the Hermite interpolation problem
\eqref{eq:basicHermiteMnf}. Because $T_q\mcM$ is a vector space, we can utilize 
the ansatz of \eqref{eq:VectorCubicHermite} but for $\V = T_q\mcM$,
\[
 \gamma(t) = a_0(t) \Delta_p + a_1(t) \Delta_q + b_0(t) \hat{v}_p + b_1(t) \hat{v}_q.
\]
Here, $ \Delta_p = \Log_q^\mcM(p),  \Delta_q = \Log_q^\mcM(q) = 0 \in T_q\mcM$ are
the normal coordinate images of the locations $p$ and $q$.
The tangent vectors
$\hat{v}_p, \hat{v}_q\in T_q\mcM$ play the role of the velocity vectors
and must be chosen such that
\begin{eqnarray}
 \dot c(t_0) &=&\frac{d}{dt}\big\vert_{t=t_0} \Exp^\mcM_q(\gamma(t)) \stackrel{!}{=} v_p = \dot f(t_0),\label{eq:tangCond0}\\
 \dot c(t_1) &=&\frac{d}{dt}\big\vert_{t=t_1} \Exp^\mcM_q(\gamma(t)) \stackrel{!}{=} v_q = \dot f(t_1)\label{eq:tangCond1}.
\end{eqnarray}
Since the interpolating curve $c(t) = \Exp^\mcM_q(\gamma(t))$ is expressed in normal coordinates centered at $q=f(t_1)$,
condition \eqref{eq:tangCond1}
is readily fulfilled by selecting $\hat{v}_q= v_q$: According to the properties of the cubic Hermite coefficient functions $a_0(t), b_0(t),b_1(t)$, the Taylor expansion of $\gamma(t)$ around $t_1$ is $\gamma(t_1+h) = h\hat{v}_q + \mathcal{O}(h^2)$. Therefore, up to first order, $\gamma(t)$ is a ray {\em emerging from the origin} $0\in T_q\mcM$ with 
direction $\hat{v}_q\in T_q\mcM$. Hence, the directional derivative of the exponential function is
\[
 \frac{d}{dt}\big\vert_{t=t_1} \Exp^\mcM_q(\gamma(t)) = \frac{d}{dt}\big\vert_{h=0} \Exp^\mcM_q(h\hat{v}_q + \mathcal{O}(h^2))= d(\Exp^\mcM_q)_0 (\hat v_q) = \hat v_q.
\]
The latter equation holds, because $d(\Exp^\mcM_q)_0 = \id_{T_q\mcM}$, \cite[\S 3, Prop. 2.9]{DoCarmo2013riemannian}.

The condition \eqref{eq:tangCond0} is more challenging, because 
the Taylor expansion of $\gamma(t)$ around $t_0$ is $\gamma(t_0+h) = \Delta_p + h\hat{v}_p + \mathcal{O}(h^2)$
and is {\em not}  a ray emerging from the origin $0\in T_q\mcM$.
As the differential $d(\Exp^\mcM_q)_{v}$ is {\em not} the identity for $v\neq 0$,
the computation of 
$\frac{d}{dt}\big\vert_{t=t_1} \Exp^\mcM_q(\Delta_p + h\hat{v}_p + \mathcal{O}(h^2))
= d(\Exp^\mcM_q)_{\Delta_p} (\hat v_p)$ is more involved.
In fact, it is related to the Jacobi fields on a Riemannian manifold, see \cite[\S 5]{DoCarmo2013riemannian},
\cite[\S 10]{Lee1997riemannian} and the upcoming Section \ref{sec:errorCurvature}.
Yet, for our purposes, it is sufficient to determine the tangent vector $\hat{v}_p$ such that 
\begin{equation}
 \label{eq:tangCond0_b}
 v_p = d(\Exp^\mcM_q)_{\Delta_p} (\hat v_p).
\end{equation}

As long as the sample points $p = \Exp^\mcM_q(\Delta_p)$ and $q$ are not conjugate,
we can make use of the fact that $\Exp^\mcM_q$ is a local diffeomorphism around $\Delta_p$, \cite[Prop. 10.11]{Lee1997riemannian}.
Hence, under this assumption, \eqref{eq:tangCond0_b} is equivalent to
\begin{eqnarray*}
 d(\Log^\mcM_q)_p(v_p) &=& \left(d(\Log^\mcM_q)_p \circ d(\Exp^\mcM_q)_{\Delta_p}\right) (\hat v_p)\\
 &=& \frac{d}{dt}\big\vert_{t=t_0} \left(\Log^\mcM_q\circ\Exp^\mcM_q\right)(\Delta_p + h\hat{v}_p)\\
 &=& \frac{d}{dt}\big\vert_{t=t_0} \id_{T_q\mcM}(\Delta_p + h\hat{v}_p) = \hat v_p.
\end{eqnarray*}
Recall that $v_p = \dot f(t_0)$ is the given sample data.
In summary, we have proved the following theorem.
\begin{theorem}
 \label{thm:cubicHermiteMnf}
 Let $t_0<t_1\in \R$ and let $f:[t_0,t_1]\rightarrow \mcM$ be a differentiable function on a Riemannian manifold $\mcM$.
 Suppose that 
 $$f(t_0) = p,\hspace{0.1cm}f(t_1)=q\in \mcM,\quad  \dot f(t_0) = v_p \in T_p\mcM,\hspace{0.1cm} \dot f(t_1) = v_q \in T_q\mcM$$
 and assume that $p$ and $q$ are not conjugate along the geodesic $t\mapsto \Exp^\mcM_q( t\Log^\mcM_q(p))$
 that connects $p$ and $q$.
 Set $\Delta_p = \Log^\mcM_q(p)$ and
 \[
  \hat v_p =  d(\Log^\mcM_q)_p(v_p) \in T_q\mcM,
  \quad \hat v_q = v_q \in T_q\mcM.
 \]
 Then
 \begin{equation}
  \label{eq:cubicHermiteMnf}
  c:[t_0,t_1] \to \mcM, \quad t\mapsto \Exp_q^\mcM\bigl(a_0(t) \Delta_p + b_0(t) \hat{v}_p + b_1(t) \hat{v}_q\bigr)
  \end{equation}
 with the cubic Hermite coefficient functions $a_0(t), b_0(t), b_1(t)$ as in  \eqref{eq:a0}-- \eqref{eq:b1}
 is a differentiable curve that solves the Hermite interpolation problem \eqref{eq:basicHermiteMnf}.
\end{theorem}
\begin{remark}
 \label{rem:compositeHermite}
 Consider an Hermite sample set $p_i = f(t_i)\in \mcM$, $v_{p_i}\in T_{p_i}\mcM$ $i=0,\ldots,k$.
 Then, by construction and in complete analogy to the Euclidean case, the composite curve
 \begin{equation}
 \label{eq:compositeHermite}
  C:[t_0, t_k]\to \mcM \quad t\mapsto c_{[t_i,t_{i+1}]}(t) \text{ for } t\in [t_i,t_{i+1}]
 \end{equation}
 that combines the local quasi-cubic spline arcs $c = c_{[t_i,t_{i+1}]}$ 
 of Theorem \ref{thm:cubicHermiteMnf}
 is of class $C^1$ and solves the Hermite manifold interpolation problem \eqref{eq:basicHermiteMnf}.
\end{remark}

\paragraph{Practical computation of $\hat{v}_p$}
%
In cases, where an explicit formula for the Riemannian logarithm is at hand, the directional derivative
 $\hat v_p  = d(\Log^\mcM_q)_p(v_p)$ can be directly computed.
For general nonlinear manifolds $\mcM$, computing the differentials of the Riemannian exponential and logarithm is rather involved.
According to \eqref{eq:RiemannExp}, \eqref{eq:RiemannLog}, \eqref{eq:diff_map}, it holds
\begin{eqnarray}
\label{eq:dLog}
  d(\Log^\mcM_q)_p :& T_p\mcM &\to T_{\Log^\mcM_q(p)}(T_q\mcM)\cong T_q\mcM,\\
\label{eq:dExp}  
  d(\Exp^\mcM_p)_0 :& T_0(T_p\mcM)\simeq T_p\mcM &\to T_{\Exp^\mcM_p(0)}\mcM =  T_p\mcM,
\end{eqnarray}
with the usual identification a linear space with its tangent space.

In order to evaluate $d(\Log^\mcM_q)_p(v_p)$, we can take any differentiable curve $\tilde \gamma(s)\subset \mcM$
that satisfies $\tilde \gamma(0) = p$ and $\dot{\tilde \gamma}(0)=v_p$.
Then,
\begin{equation}
 \label{eq:vhat_comp}
  d(\Log^\mcM_q)_p(v_p) = d(\Log^\mcM_q)_{\tilde{\gamma}(0)}(\dot{\tilde \gamma}(0))
  = \frac{d}{ds}\big\vert_{s=0} \Log^\mcM_q\left(\tilde{\gamma}(s)\right).
\end{equation}
An obvious choice is $\tilde \gamma(s) = \Exp_p(sv_p)\subset \mcM$.
The final equation for computing $\hat v_p$ as required by \eqref{eq:tangCond0} is
\begin{equation}
 \label{eq:vhat_final}
 \hat v_p = \frac{d}{ds}\big\vert_{s=0}\left(\Log^\mcM_q\circ \Exp^\mcM_p\right)(sv_p).
 \footnote{\text{Due to the different base points, 
 this composition of $\Log^\mcM_q$ and $\Exp^\mcM_p$ is not the identity.}}
\end{equation}
The composite map $\Log^\mcM_q\circ \Exp^\mcM_p:T_p\mcM\supset D_0 \to T_q\mcM$ is in fact a transition function 
for the normal coordinate charts. It is defined on an open subset of
a Hilbert space and  maps to a Hilbert space, see \cite[Fig. 1.6, p. 12]{Lee2012smooth} for an illustration.
Hence, we can approximate the directional derivative
$ \hat v_p =  \frac{d}{ds}\big\vert_{s=0}\left(\Log^\mcM_q\circ \Exp^\mcM_p\right)(sv_p)$
via finite difference approaches:
\begin{equation}
 \label{eq:FD_dirdiff}
  \hat v_p 
 =
 \frac{\left(\Log^\mcM_q\circ \Exp^\mcM_p\right)(h v_p) - \left(\Log^\mcM_q\circ \Exp^\mcM_p\right)(-h v_p)}{2h}
 + \mathcal{O}(h^2).
\end{equation}
\subsection{Computational effort and preliminary comparison to other methods}
\label{sec:comparison}
Computationally, the most involved numerical operations are the evaluations
of Riemannian $\Log$- and $\Exp$-mappings.
Therefore, as in \cite{GouseMassartAbsil2018}, we measure the computational effort associated with the Hermite interpolation method as the number of such function evaluations.

{\em Constructing} a quasi-cubic Hermite interpolant as in Remark \ref{rem:compositeHermite} requires
on each subinterval $[t_i,t_{i+1}]$
\begin{itemize}
 \item one Riemannian logarithm to compute $\Delta_p = \Log^\mcM_q(p)$,
 \item two Riemannian $\Log$- and $\Exp$-evaluations for the central difference approximation of \eqref{eq:FD_dirdiff},
\end{itemize}
which results in a total of $3k$ Riemannian $\Log$-evaluations 
and $2k$ Riemannian $\Exp$-evaluations for the whole composite curve.
The data to represent the curve \eqref{eq:compositeHermite} can be precomputed and stored.

{\em Evaluating} a quasi-cubic Hermite interpolant
at time $t$ requires {\em a single} Riemannian $\Exp$-evaluation.

As mentioned in the introduction, B{\'e}zier-like approaches may be used to tackle the 
the Hermite interpolation problem \eqref{eq:basicHermiteMnf}.
This requires a cubic degree and at least four control points on each sub-interval $[t_i,t_{i+1}]$
to impose the derivative constraints,
see Fig. \ref{fig:CubicBezier}.
The most efficient of such methods in \cite{GouseMassartAbsil2018} 
requires $\mathcal{O}(k^2)$ Riemannian $\Log$-evaluations for constructing the curve data.
Evaluating the curve at time $t$ requires
 $3$ Riemannian $\Exp$-evaluations plus $1$ Riemannian $\Log$-evaluation \cite[Prop. 5.10]{GouseMassartAbsil2018}.
\section{Error propagation}
\label{sec:errorCurvature}
The approach introduced in Section \ref{sec:HermiteInterp2Mnf} follows the standard principle of 
(1) mapping the sampled data onto the tangent space, (2) performing data processing (in this case, interpolation) in the tangent space, (3) mapping the result back to the curved manifold.
In this section, we perform a general qualitative analysis of the behavior of the actual errors on the manifold
in question in relation to the data processing errors that accumulate in the tangent space.
In particular, this allows to obtain error estimates for any manifold interpolation procedure based on the above standard principle
and also applies to other data processing operations that subordinate to this pattern.
In essence, the error propagation is related to the manifold's curvature via a standard result from differential geometry
on the spreading of geodesics \cite[Chapter 5, \S 2]{DoCarmo2013riemannian}.
%
\begin{theorem}
 \label{thm:ErrorMain}
 Let $\mcM$ be a Riemannian manifold, let $q\in \mcM$ and consider tangent vectors
 $\Delta, \tilde{\Delta}\in T_q\mcM$, which are to be
 interpreted as exact datum and associated approximation.
 Write $\delta = \|\Delta\|$, $\tilde{\delta} = \|\tilde{\Delta}\|$,
 where it is understood that the norm is that of $T_q\mcM$.
 Assume that $\delta, \tilde{\delta} < 1$.
 Let $\sigma=\text{span}(\Delta, \tilde{\Delta})\subset T_q\mcM$ and 
 let $K_q(\sigma)$ be the sectional curvature at $q$ with respect to the $2$-plane $\sigma$.
 
 If $s_0 = \angle(\tilde{\Delta},\Delta)$ is the angle between $\tilde{\Delta}$ and $\Delta$,
 then the Riemannian distance between the manifold locations $p = \Exp_q^\mcM(\Delta)$ and
 $\tilde{p} = \Exp_q^\mcM(\tilde{\Delta})$ is
 \begin{equation}
 \label{eq:ErrordistFinal}
 \dist_\mcM(p, \tilde{p}) \leq |\delta-\tilde{\delta}| + s_0\delta 
   (1 - \frac{K_q(\sigma)}{6}\delta^2 + o(\delta^2)) + \mathcal{O}(s_0^2),
 \end{equation}
 with the underlying assumption that all data is within the injectivity radius at $q$. 
\end{theorem}
%
%
\begin{figure}[ht]
\centering
\includegraphics[width=0.7\textwidth]{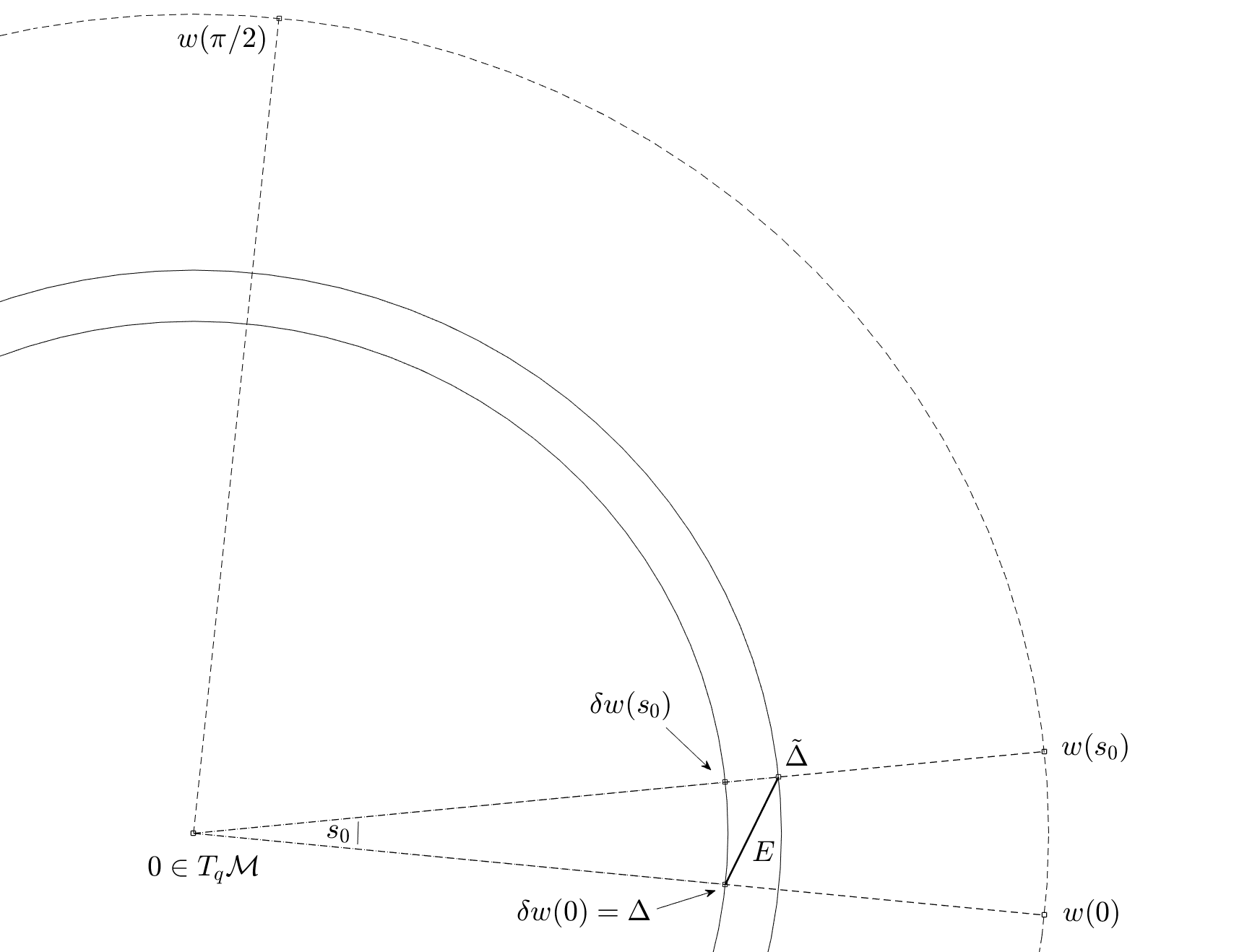}
\caption{A
}
\label{fig:PlotErrEst}
\end{figure}
%
\begin{proof}
 Formally, it holds  $\dist_\mcM(p, \tilde{p}) = \|\Log_p^\mcM(\tilde{p})\|_{T_p\mcM}$.
 However, the data is given in normal coordinates around $q\in\mcM$ and not around $p$ (nor $\tilde{p}$).
 The Riemannian exponential is a radial isometry (lengths of rays starting from the origin of the tangent space
 equal the lengths of the corresponding geodesics).
 Yet, it is not an isometry
 so that  $\dist_\mcM(p, \tilde{p})\neq \|E\|$, unless $\mcM$ is flat.
 Therefore, we will estimate the distance $\dist_\mcM(p, \tilde{p})$
 against a component that corresponds to the length of a ray in $T_q\mcM$
 and a circular segment in $T_q\mcM$.
 
 To this end, introduce an orthonormal basis ${w, w^\bot}$ for the plane 
 $\sigma=\text{span}(\Delta, \tilde{\Delta})\subset T_q\mcM$ via
 \[
  w:=\frac{\Delta}{\|\Delta\|}, \quad w^\bot = \frac{\tilde{\Delta} - \langle w, \tilde{\Delta}\rangle w}{\|\tilde{\Delta} - \langle w, \tilde{\Delta}\rangle w\|}.
 \]
 The circular segment $w(s)$ of the unit circle in the $\sigma$-plane that starts from $w=w(0)$ and ends in $w^\bot=w(\pi/2)$ can be parameterized via the curve
 \[
  w:[0,\pi/2]\to T_q\mcM, \quad s \mapsto w(s) = \cos(s)w + \sin(s)w^\bot.
 \]
 Let $s_0\in [0,\pi/2]$ be the angle such that $\tilde{\Delta} = \tilde{\delta}w(s_0)$. This setup is illustrated in Figure \ref{fig:PlotErrEst},
 where the outer dashed circular arc indicates the unit circle and the solid circular arcs are the circles of radius $\delta$ and $\tilde{\delta}$, respectively.
 By the triangle inequality,
 \begin{eqnarray}
   \label{eq:triangle0}
  \dist_\mcM(p, \tilde{p})&=& \dist_\mcM\left(\Exp_q^\mcM(\delta w(0)), \Exp_q^\mcM( \tilde{\delta} w(s_0))\right)\\
  \label{eq:triangle1}
  &\leq& \dist_\mcM\left(\Exp_q^\mcM(\delta w(0)), \Exp_q^\mcM( \delta w(s_0))\right)\\
   \label{eq:triangle2}
  &&+ \dist_\mcM\left(\Exp_q^\mcM(\delta w(s_0)), \Exp_q^\mcM( \tilde{\delta} w(s_0))\right).
 \end{eqnarray}
 Since the points $\delta w(s_0)$ and $\tilde{\delta} w(s_0)=\tilde\Delta$ are on a ray
 that emerges from the origin in $T_q\mcM$, the distance term in line \eqref{eq:triangle2}
 is exactly $|\delta - \tilde{\delta}|$, see Figure \ref{fig:PlotErrEst}.
 Note that $\Exp_q^\mcM(\delta w(0))=\Exp_q^\mcM(\Delta)=p$.
 Hence, the distance term in line \eqref{eq:triangle1} is
 \begin{eqnarray*}
  \dist_\mcM\left(p, \Exp_q^\mcM( \delta w(s_0))\right)
  &=& \|\Log_p^\mcM( \Exp_q^\mcM( \delta w(s_0)))\|_{T_p\mcM}.
 \end{eqnarray*}
 A Taylor expansion centered at $s=0$ of the transition function along the circular segment
 $s\mapsto \left(\Log_p^\mcM\circ \Exp_q^\mcM\right)( \delta w(s)))$
 gives
  \begin{eqnarray*}
  \Log_p^\mcM( \Exp_q^\mcM( \delta w(s_0)))
  &=& \Log_p^\mcM( \Exp_q^\mcM( \delta w(0))) \\
  && + s_0\frac{d}{ds}\big\vert_{s=0} \left(\Log_p^\mcM\circ \Exp_q^\mcM\right)( \delta w(s)))
  +\mathcal{O}(s_0^2)\\
  &=&  \Log_p^\mcM(p) + s_0 d(\Log_p^\mcM)_p \circ d(\Exp_q^\mcM)_\Delta (\delta \dot w(0)) +\mathcal{O}(s_0^2)\\
  &=& 0 + s_0 d(\Exp_q^\mcM)_{\delta w} (\delta w^\bot) +\mathcal{O}(s_0^2).
 \end{eqnarray*}
To arrive at the last line, $ d(\Log_p^\mcM)_p = \id_{T_p\mcM}$ was used, which follows from
the standard result  $d(\Exp_p^\mcM)_0 = \id_{T_p\mcM}$ \cite[\S 3, Prop. 2.9, p. 65]{DoCarmo2013riemannian}
with the usual identification of $T_p\mcM \cong T_p(T_p\mcM)$, cf. \eqref{eq:dLog}, \eqref{eq:dExp}.

By construction, $\delta \mapsto d(\Exp_q^\mcM)_{\delta w} (\delta w^\bot)$ is a Jacobi field along
the geodesic ray that starts from $q = \Exp_q^\mcM(0)$ with unit velocity $w\in T_q\mcM$,
see \cite[\S 5, Prop. 2.7, p. 114]{DoCarmo2013riemannian}.
Moreover, $w=(0), w^\bot=w(\pi/2)$ constitute an orthonormal basis of the $\sigma$-plane in $T_q\mcM$.
Therefore, the results \cite[\S 5, Cor. 2.9, Cor. 2.10, p. 115]{DoCarmo2013riemannian} apply
and give
\[
 \|d(\Exp_q^\mcM)_{\delta w} (\delta w^\bot)\| = \delta - \frac{K_q(\sigma)}{6}\delta^3 + o(\delta^3).
\]
In summary,
\begin{equation*}
  \dist_\mcM(p, \tilde{p}) =|\delta - \tilde{\delta}|+ 
 \delta s_0\left(1 - \frac{K_q(\sigma)}{6}\delta^2 + o(\delta^2)\right) +\mathcal{O}(s_0^2),
\end{equation*}
which establishes the theorem.
\end{proof}
%
\begin{remark}
\begin{enumerate}[(i)]
 \item The the approximation error in the tangent space $\epsilon:=\|E\| = \|\tilde{\Delta}-\Delta \|$
 can be related by elementary trigonometry to the angle $s_0 = \angle(\tilde{\Delta},\Delta)$.
 It holds $\epsilon\geq \delta\|w(0)-w(s_0)\|$, see Fig. \ref{fig:PlotErrEst}.
 Moreover, $s_0 = 2\arcsin\left(\frac{\|w(0)-w(s_0)\|}{2}\right)$.
 Thus,
 $$\delta s_0 \leq  2\delta\arcsin\left(\frac{\epsilon}{2\delta}\right) = \epsilon +
 \mathcal{O}\left(\frac{\epsilon^3}{(2\delta)^2}\right).$$
 In regards of practical applications, it is safe to assume $\epsilon < \delta$. Then, 
 in terms of the error $\epsilon$, the distance estimate \eqref{eq:ErrordistFinal} reads
 \begin{equation}
  \label{eq:ErrordistFinal_err}
  \dist_\mcM(p, \tilde{p}) =|\delta - \tilde{\delta}| +  
  \epsilon\left(1 - \frac{K_q(\sigma)}{6}\delta^2 + o(\delta^2)\right) +\mathcal{O}(\epsilon^2).
 \end{equation}
 \item
 If we travel from $\Delta$ to $\tilde{\Delta}$ in the tangent space on the corresponding curves
 as in the proof of Theorem \ref{thm:ErrorMain},
 i.e. first along the circular arc from $\Delta =\delta w(0)$ to $\delta w(s_0)$
 and than along the ray from $\delta w(s_0)$to $\tilde{\delta} w(s_0)$,
 then we cover precisely a distance of 
 $|\delta -\tilde{\delta}| + \delta s_0$.
 Comparing this with \eqref{eq:ErrordistFinal}, we see that the corresponding distances of the 
 manifold images are (asymptotically) 
 $\left\{\begin{array}{l}
   \text{shorter,}\\
   \text{longer,}
  \end{array}
  \right\}
 $
 if $\mcM$ features 
  $\left\{\begin{array}{l}
   \text{positive}\\
   \text{negative}\\
  \end{array}
  \right\}
 $ sectional curvatures.
 The underlying principle is the well-known effect that geodesics on positively curved spaces spread apart less than
 straight rays in a flat space, while they spread apart more on negatively curved spaces, see \cite[\S 5, Remark 2.11, p. 115/116]{DoCarmo2013riemannian}.
 
 From the numerical point of view, this means that data processing operations that work in the tangent space followed by a 
 transition to the manifold are rather well-behaved on manifolds of positive curvature, while the opposite holds on negatively curved manifolds.
 In Section \ref{sec:Stiefel_Hermite}, we will show an illustration of Theorem \ref{thm:ErrorMain} 
 on an interpolation problem on the compact Stiefel manifold.
 \end{enumerate}
\end{remark}
%
With the help of Theorem \ref{thm:ErrorMain}, explicit error bounds for manifold interpolation methods
can be obtained. For example, cubic Hermite interpolation comes with a standard error bound \cite[Thm 7.16]{hohmann2003numerical}
that applies to the interpolant in the tangent space.
This can be forwarded to a manifold error via Theorem \ref{thm:ErrorMain}.
%
%
%

%
%
%
\section{Cubic Hermite interpolation of column-orthogonal matrices}
\label{sec:Stiefel_Hermite}
The set of column-orthogonal matrices
\[
  St(n,r):= \{U \in \R^{n\times r}| \quad U^TU = I_r\}
\]
is the compact homogeneous matrix manifold known as the (compact) {\em Stiefel manifold}.
This section reviews the essential aspects of the numerical treatment of
Stiefel manifolds. For more details, see \cite{AbsilMahonySepulchre2008, EdelmanAriasSmith1999,Zimmermann_MORHB2019}.

The {\em tangent space} $T_USt(n,r)$ at a point $U \in St(n,r)$
can be thought of as the space of velocity vectors of differentiable curves on $St(n,r)$
passing through $U$:
\[
  T_USt(n,r)=\{\dot{c}(t_0)| c:(t_0-\epsilon, t_0+\epsilon)\rightarrow St(n,r), c(t_0)=U\}.
\]
For any matrix representative $U\in St(n,r)$,
the tangent space of $St(n,r)$ at $U$ is
\[
  T_USt(n,r) = \left\{\Delta \in \R^{n\times r}|\quad U^T\Delta = -\Delta^TU\right\}\subset \R^{n\times r}.
\]
Every tangent vector $\Delta \in T_USt(n,r)$  may be written as
\begin{equation}
\label{eq:tang2}
  \Delta = UA + (I-UU^T)T, \quad A \in \R^{r\times r} \mbox{ skew}, \quad T\in\R^{n\times r} \mbox{ arbitrary}.
\end{equation}
%
%
The dimension of both $T_USt(n,r)$ and $St(n,r)$ is $nr -\frac{1}{2}r(r+1)$.

Each tangent space carries an inner product 
$\langle \Delta, \tilde{\Delta}\rangle_U = tr\left(\Delta^T(I-\frac{1}{2}UU^T)\tilde{\Delta}\right)$
\tc{with corresponding norm $\| \Delta\|_U = \sqrt{\langle \Delta, \Delta\rangle_U}$.}
This is called the {\em canonical metric} on $T_USt(n,r)$.
It is derived from the 
quotient space representation $St(n,r) = O(n)/O(n-r)$
that identifies two square orthogonal matrices in $O(n)$
as the same point on $St(n,r)$, if their first $r$ columns coincide \cite[\S 2.4]{EdelmanAriasSmith1999}.
For a condensed introduction to quotient spaces, see \cite[\S 2.5]{Zimmermann_MORHB2019}.
Endowing each tangent space with this metric (that varies differentiably in $U$)
turns $St(n,r)$ into a {\em Riemannian manifold}.
The associated sectional curvature is non-negative and is bounded by $0\leq K_U(\sigma) \leq \frac{5}{4}$
for all $U\in St(n,r)$ and all two-plans $\sigma = \text{span}(\Delta, \tilde{\Delta})\subset T_USt(n,r)$, \cite[\S 5]{Rentmeesters2013}.

Given a start point $U\in St(n,r)$ and an initial velocity 
$\Delta\in T_USt(n,r)$ the Stiefel geodesic $c_{U, \Delta}$ (and thus the Riemannian exponential)
is
\begin{equation}
 \label{eq:Stexp}
 c_{U, \Delta}(t) = Exp_{U}^{St}(t\Delta) 
 = (U,Q)\exp_m\left(t\begin{pmatrix}A & -R^T\\ R & 0\end{pmatrix}\right)\begin{pmatrix}I_r\\ 0\end{pmatrix},
\end{equation}
where
\[
   \Delta = UU^T\Delta + (I-UU^T)\Delta \stackrel{(\mbox{\footnotesize QR-decomp. of } (I-UU^T)\Delta)}{=} UA + QR
\]
is the decomposition of the tangent velocity into its horizontal and vertical component
with respect to the base point $U$, \cite{EdelmanAriasSmith1999}.
Because $\Delta$ is tangent, $A=U^T\Delta\in\R^{r\times r}$ is skew.
The Riemannian Stiefel logarithm can be computed with the algorithm of \cite{StiefelLog_Zimmermann2017}.

\subsection{Differentiating the Stiefel exponential}
In this section, we compute the directional derivative of the Stiefel exponential 
\begin{equation}
 \label{eq:diffStiefel}
  \frac{d}{dt}\big\vert_{t=0} \Exp_U^{St}(\Delta_0 + tV), \quad \Delta_0, V\in T_USt(n,r).
\end{equation}
This is important for two reasons.
\begin{enumerate}
 \item {\em Differentiable gluing of interpolation curves.}
 Consider a manifold data set $t_i, p_i = f(t_i)$, $i=0,\ldots, j, j+1,\ldots, k$,
 where the Riemannian distance, say, of the sample points $p_j$ and $p_{0}$ 
 and $p_j$ and $p_{k}$ exceeds the  injectivity radius of $\mcM$ at $p_j$.
 Then, simple tangent space interpolation with mapping the data set to $T_{p_j}\mcM$
 is not possible.
 A remedy is to split the data set at $p_j$ and to compute two interpolation curves,
 one for the sample set  $t_i, p_i = f(t_i)$, $i=0,\ldots, j$
 and one for the sample set $t_i, p_i = f(t_i)$, $i=j, j+1,\ldots, k$.
 With the canonical method of tangent space interpolation, the curves have the expressions
 $c_1(t) = \Exp_{p_{\lfloor j/2 \rfloor}}(\sum_{i=0}^ja_i(t) \Delta_i)$ and
 $c_2(t) = \Exp_{p_{j+\lfloor j/2 \rfloor}}(\sum_{i=j}^ka_i(t) \Delta_i)$,
 where $\Delta_i = \Log^{St}_{p_{\lfloor j/2 \rfloor}}(p_i)$ for $c_1$ and
 $\Delta_i = \Log^{\mcM}_{p_{j + \lfloor j/2 \rfloor}}(p_i)$ for $c_2$.
 Concatenating the curves $c_1$, $c_2$ will result in a non-differentiable kink 
 at the intersection location $p_j$, where $c_1$ ends and $c_2$ starts.
 In order to avoid this, one can compute the derivative $\dot c_1(t_j)$
 and use $\dot c_1(t_j) = \dot c_2(t_j)$ as an Hermitian derivative sample when constructing $c_2$.
 For obtaining $\dot c_1(t_j)$, a derivative of the form of \eqref{eq:diffStiefel} must be computed.
 \item {\em Method validation.}
 The cubic Hermite manifold interpolation method of Theorem \ref{thm:cubicHermiteMnf} requires the computation of
 $\hat v_p  = d(\Log^\mcM_q)_p(v_p)$.
 As was mentioned in Section \ref{sec:HermiteInterp2Mnf}, the differential of the $\Log$-mapping cannot be computed explicitly
 for general manifolds $\mcM$.
 In order to assess the numerical quality of a finite-differences approximation,
 we can first compute $\hat v_p$ by \eqref{eq:FD_dirdiff} and then recompute \eqref{eq:tangCond0_b}
 $$v_{p,rec} = d(\Exp^\mcM_q)_{\Delta_p} (\hat v_p) = \frac{d}{dt}\big\vert_{t=0}\Exp^\mcM_q(\Delta_p+ t\hat v_p).$$
 The numerical accuracy is assessed via the error
 \begin{equation}
  \label{eq:error_hatv}
  \frac{\|v_{p,rec} - v_p\|_p}{\|v_p\|_p}.
 \end{equation}
 Again, a derivative of the form of \eqref{eq:diffStiefel} must be computed.
\end{enumerate}
Now, let us address the derivative \eqref{eq:diffStiefel} for $\Delta_0, V\in T_USt(n,r)$.
The underlying computational obstacle is that the exponential law does not hold for the matrix exponential 
and two non-commuting matrices $\exp_m(X+tY)\neq \exp_m(X)\exp_m(tY)$.
Write $\Delta(t) = \Delta_0 + tV$ and let 
$Q(t)R(t) = (I-UU^T)\Delta(t)$ be the $t$-dependent QR-decomposition of the tangent space curve.
Moreover, $A(t) := U^T\Delta(t)$ and $\dot A(0) =  U^TV$. Then, by the product rule,
\begin{eqnarray*}
 \frac{d}{dt}\big\vert_{t=0} \Exp_U^{St}(\Delta(t))
 &=&\frac{d}{dt}\big\vert_{t=0}
 (U,Q(t))\exp_m\left(\begin{pmatrix}A(t) & -R^T(t)\\ R(t) & 0\end{pmatrix}\right)
                    \begin{pmatrix}I_r\\ 0\end{pmatrix}\\
 &=&(0, \dot Q(0)) \exp_m\left(\begin{pmatrix}A(0) & -R^T(0)\\ R(0) & 0\end{pmatrix}\right)\begin{pmatrix}I_r\\ 0\end{pmatrix}\\
 && + (U,Q(0))\frac{d}{dt}\big\vert_{t=0} \exp_m\left(\begin{pmatrix}A(t) & -R^T(t)\\ R(t) & 0\end{pmatrix}\right)\begin{pmatrix}I_r\\ 0\end{pmatrix}.
\end{eqnarray*}

Introduce the matrix function
$M(t) = \begin{pmatrix}A(t) & -R^T(t)\\ R(t) & 0\end{pmatrix}$.
%
It is sufficient to compute $d(\exp_m)_{M(0)}(\dot M(0)) = \frac{d}{dt}\big\vert_{t=0} \exp_m(M(0) + t\dot M(0))$.\footnote{This is a common problem in 
Lie group theory, see \cite[\S 5.4]{Hall_Lie2015}. The solution is formally an infinite sequence of nested commutator products in $[M,\dot M] =M\dot M -\dot M M$, 
$$\frac{d}{dt}\big\vert_{t=0} \exp_m(M + t\dot M) = \exp_m(M)\left(\dot M- \frac{1}{2!}[M,\dot M] + \frac{1}{3!}[M,[M,\dot M]]-\cdots \right)$$
}
In the following, we often omit the parameter $t$ with the implicit understanding that all quantities are evaluated at $t=0$.
By Mathias' Theorem \cite[Thm 3.6, p. 58]{Higham:2008:FM}, it holds
\begin{equation}
 \label{eq:Mathias}
 \exp_m \begin{pmatrix}[c|c] M & \dot M\\\hline  0 & M\end{pmatrix}
 =  \begin{pmatrix}[c|c] \exp_m(M) & \frac{d}{dt}\big\vert_{t=0} \exp_m(M + t\dot M)\\\hline 0 & \exp_m(M)\end{pmatrix}.
\end{equation}
Hence, for data stemming from $St(n,r)$, a $(4r\times 4r)$-matrix exponential must be computed.
However, the advantage is that $\exp_m(M)$ and $\frac{d}{dt}\big\vert_{t=0} \exp_m(M + t\dot M)$ are obtained in one go and both are needed for evaluating \eqref{eq:diffStiefel}. Moreover, usually $r\ll n$ in practical applications.
For details and alternative algorithms for computing  $\frac{d}{dt}\big\vert_{t=0} \exp_m(M + t\dot M)$,
see \cite[\S 10.6]{Higham:2008:FM}.
In summary:
\begin{lemma}
\label{lem:Stdiff}
 With all quantities as introduced above,
 let
 \[
 \exp_m \begin{pmatrix}[c|c] M & \dot M\\\hline  0 & M\end{pmatrix} = 
  \begin{pmatrix}[c|c]
   \begin{pmatrix}
   E_{11} & E_{12}\\
   E_{21} & E_{22}
   \end{pmatrix}
   &
      \begin{pmatrix}
   D_{11} & D_{12}\\
   D_{21} & D_{22}
   \end{pmatrix}\\ \hline
   \mathbf{0}
   &
      \begin{pmatrix}
   E_{11} & E_{12}\\
   E_{21} & E_{22}
   \end{pmatrix}
  \end{pmatrix}
 \]
 be written in terms of subblocks of size $r\times r$.
 Then
 \begin{equation}
 \label{eq:Stdiff_final}
  \frac{d}{dt}\big\vert_{t=0} \Exp_U^{St}(\Delta(t))
  = \dot Q E_{21} + UD_{11} + QD_{21}.
 \end{equation}
\end{lemma}
The derivatives of the QR-factors of 
the decomposition $Q(t)R(t) = (I-UU^T)\Delta(t)$ that are required to compute $\dot Q$ and $\dot M = \begin{pmatrix}\dot A & -\dot R^T\\ \dot R & 0\end{pmatrix}$ can be obtained from Alg.  \ref{alg:QRdiff}.
\subsection{Alternative options for Hermite data preprocessing on St(n,r)}
As outlined in Section \ref{sec:HermiteInterpMnf}, the Hermite interpolation problem
with local Stiefel sample data 
$f(t_0)= U, \dot f(t_0) = \Delta$, $f(t_1)= \tilde{U}, \dot f(t_1) = \tilde{\Delta}$,
 requires us to translate the derivative samples to a common tangent space.
On $St(n,r)$, this amounts to compute 
\begin{equation*}
  \hat{\Delta} = \frac{d}{ds}\big\vert_{s=0} \Log^{St}_{\tilde{U}}\left(\tilde{\gamma}(s)\right)
\end{equation*}
for some differentiable curve $\tilde \gamma(s)\subset St(n,r)$
that satisfies $\tilde \gamma(0) = U$ and $\dot{\tilde \gamma}(0)=\Delta$.
There are other option than $\tilde\gamma(s) = \Exp_{U}^{St}(s\Delta)$ with this property,
which might be cheaper to evaluate, depending on the context:
For a skew-symmetric $M_0$,
the {\em Cayley transformation} \cite[eq. (7)]{WenYin2013}, \cite[p. 284]{Bhatia1997},
\begin{equation}
 \label{eq:cayley}
 M:s \mapsto M(s) = \left(I + \frac{s}{2}M_0\right)\left(I - \frac{s}{2}M_0\right)^{-1} = I + sM_0 + \frac{s^2}{2}M^2_0 + \dots
\end{equation}
produces a curve of orthogonal matrices that 
matches the matrix exponential $s\mapsto \exp_m(s M_0)$
up terms of order $\mathcal{O}(s^2)$.
As a consequence,
\[
 M_0=\frac{d}{ds}\big\vert¦_{s=0} \exp_m(s M_0)
 = \frac{d}{ds}\big\vert¦_{s=0} M(s)
\]
and the matrix curve $s\mapsto (U,Q) M(s)\begin{pmatrix}I_r\\ 0\end{pmatrix}$ 
with $M(s)$ from \eqref{eq:cayley} based on $M_0 = \begin{pmatrix}A & -R^T\\ R & 0\end{pmatrix}$ as in \eqref{eq:Stexp}
may be used as the curve $\tilde{\gamma}(s)$ in \eqref{eq:vhat_comp}, \eqref{eq:vhat_final}
instead of the Stiefel exponential \eqref{eq:Stexp}.

Another option is to use {\em retractions} as a replacement for the Riemannian exponential \cite[\S 4.1]{AbsilMahonySepulchre2008}.
By definition, the differential of a retraction map at the origin of the tangent space is the identity map and thus coincides with  the differential of the Riemannian exponential at the origin.
Suitable matrix curves $s\mapsto \tilde{\gamma}(s)$ that match $\Exp_{U}^{St}(s\Delta)$ up to terms of first order based on 
Stiefel retractions are
\begin{eqnarray*}
 &\tilde{\gamma}: s\mapsto (U + s\Delta)\Phi (I+s^2\Lambda)^{-\frac{1}{2}}\Phi^T, & \quad 
 \Phi \Lambda\Phi^T \stackrel{\text{\small EVD}}{=} \Delta^T\Delta \\
 &\tilde{\gamma}: s\mapsto qr(U + s\Delta), &\quad  \text{\small (compact qr-decomposition)}
\end{eqnarray*}
see \cite[Example 4.1.3]{AbsilMahonySepulchre2008}.\\
\textbf{A word of caution:}
With the QR-based retraction, there is the challenge of 
computing a {\em differentiable} QR-path. Numerical QR-algorithms in high-level programming environments like MATLAB or SciPy
might provide discontinuous matrix paths, e.g., because of different internal pivoting strategies.
%
%
%
%
\section{Examples and experimental results}
\label{sec:experiments}
%
%
%
%
%
%
%
%
%
%
%
In this section, we conduct various numerical experiments that put the theoretical findings in perspective.
All examples are coded and performed in the SciPy programming environment \cite{SciPy2001}.
\subsection{The numerical accuracy of the derivative translates}
\label{sec:FD_errors}
Before we start with the actual interpolation problems, we assess the numerical accuracy of
the process of mapping a velocity sample $v_U\in T_USt(n,r)$ to a tangent velocity $\hat{v}_U \in T_{\tilde{U}}St(n,r)$
for two different Stiefel locations $U, \tilde{U}\in St(n,r)$.
This requires the numerical computation of
$\hat v_U  = d(\Log^{St}_{\tilde{U}})_U(v_U)$ with the help of central finite differences as in \eqref{eq:FD_dirdiff}.
 
Then, we reverse this process to recover the original input via
$$v_{U,rec} =  \frac{d}{dt}\big\vert_{t=0}\Exp^{St}_{\tilde{U}}( \Log_{\tilde{U}}^{St}(U) + t\hat v_U).$$
To this end, we utilize formula \eqref{eq:Stdiff_final} of Lemma \ref{lem:Stdiff}.
Then we compute the error
$
 \mathcal{E} = \frac{\|v_{U,rec} - v_U\|_F}{\|v_U\|_F}.
$
As data points, we use samples of the Stiefel function 
$\mu\mapsto U(\mu)\in St(n,r)$, $n=1001, r=6$ that is featured in the upcoming Section \ref{sec:CHSVD_fun}:
More precisely, $U = U(0.9)$, $\tilde{U} = U(1.4)$.
The tangent direction to be translated is chosen as $v_U = Log_{U}^{St}(U(1.9))\in T_USt(n,r)$.
All Riemannian log computations are performed with the Algorithm of \cite{StiefelLog_Zimmermann2017}
and a numerical convergence threshold of $\tau = 10^{-14}$.
The next table shows the reconstruction error $\mathcal{E}$ versus the finite difference step size $h$ used in \eqref{eq:FD_dirdiff}.
\begin{center}
\begin{tabular}{l|l|l|l|l|l|l|l}
 step size     $h$   &  $10^{-2}$ & $10^{-3}$ & $10^{-4}$ & $10^{-5}$ & $10^{-6}$ & $10^{-7}$\\  \hline
 error $\mathcal{E}$ &    1.2e-8  &   1.2e-10 &   4.3e-12 &   4.2e-11 &   4.1e-10 & 5.0e-9
\end{tabular}.
\end{center}
Even though there are various numerical processes involved (matrix exp, matrix log, numerical QR-differentiation, iterative Stiefel logarithm etc.)
the accuracy of the finite difference approach is surprisingly high.
In the following experiments, a step size of $h=10^{-4}$ is used to calculate \eqref{eq:FD_dirdiff}.
\subsection{Hermite interpolation of the Q-factor of a QR-decomposition}
\label{sec:CHQR}
As a first example, consider a cubic matrix polynomial
\[ 
 Y(t) = Y_0 + t Y_1 + t^2Y_2 + t^3Y_3, \quad Y_i \in \R^{n\times r}, n=500, r=10.
\]
The matrices $Y_i$ were produced as random matrices with entries uniformly sampled from $[0,1]$ for $Y_0$,
entries uniformly sampled from $[0,0.5]$ for $Y_1,Y_2$ and from $[0,0.2]$ for $Y_3$.
The $t$-dependent $QR$-decomposition is
 \begin{equation*}
  Y(t) = Q(t)R(t),\quad 
  \dot Y(t)= \dot Q(t)R(t) + Q(t)\dot R(t).
 \end{equation*}
The matrix curve $Y(t)$ is sampled at $6$ Chebychev roots in $[-1.1, 1.1]$.\footnote{ The Chebychev locations read $t_0\approx -1.0625$, $t_1\approx -0.7778$, $t_2 \approx -0.2847$, $t_3\approx0.2847$, $t_4 \approx  0.7778$, $t_5\approx 1.0625$.}
At each sample point $t_i$ the Q-factor $Q(t_i)$ of the QR-decomposition computed. 
The corresponding derivative $\dot Q(t_i)$ is obtained from Alg. \ref{alg:QRdiff} in Appendix \ref{app:diffQR}.
This constitutes the Hermite sample data set
\[
 Q(t_i)\in St(n,r), \quad \dot Q(t_i)\in T_{Q(t_i)}St(n,r), \quad i=0,\ldots,5.
\]
For comparison, the following interpolation schemes are conducted.
 \begin{itemize}
  \item {\em Quasi-linear interpolation:} In this case, the Stiefel samples are connected by geodesics as described in \cite[\S 3.1]{Zimmermann_MORHB2019}. No derivative information is used. This is the manifold version of linear interpolation.
  \item {\em Tangent space interpolation:} In this case, all data is mapped to single tangent space attached at 
  $Q(t_j)$, $j=\lfloor k/2 \rfloor$, where $k$ is the number of sample points.
  Then, RBF interpolation is performed on tangent vectors as described in \cite{AmsallemFarhat2008}, \cite[\S 3.1]{Zimmermann_MORHB2019}.
  As an RBF, the inverse multiquadric is selected. No derivative information is used.
  \item {\em Hermite quasi-cubic interpolation}, as introduced in Section \ref{sec:HermiteInterp2Mnf}.
 \end{itemize}
 %
%
\begin{figure}[ht]
\centering
\includegraphics[width=1.0\textwidth]{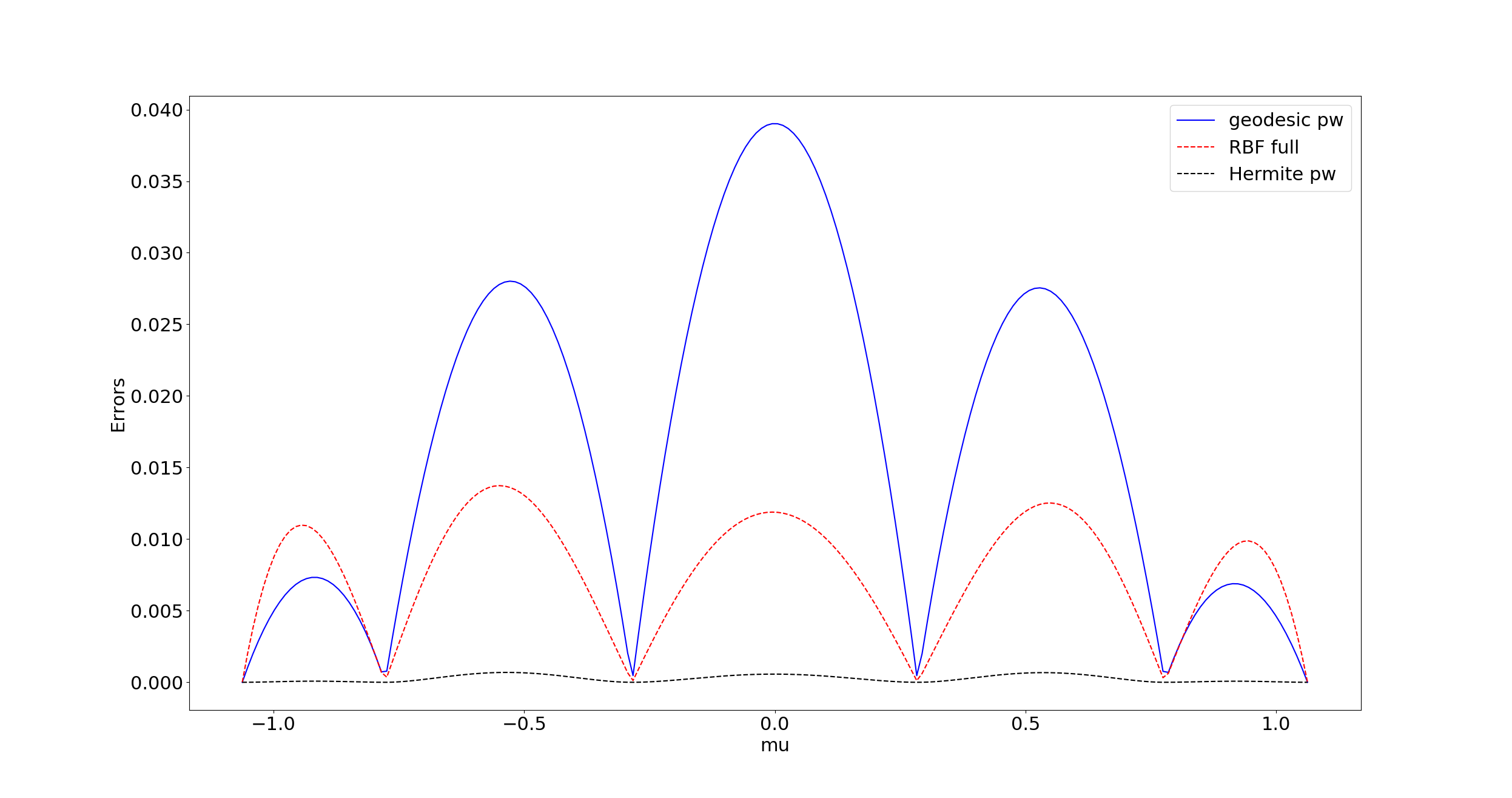}
\label{fig:Errors_analytic3_interp_n500p10_Cheby6_-1p1_1p1_invquad}
\caption{Relative errors of the various interpolation approaches in the matrix Frobenius norm
for the experiment of Section \ref{sec:CHQR}.}
\end{figure}
%
Since the quasi-linear and the quasi-cubic approach rely on piece-wise splines, it is only the `global' tangent space interpolation that benefits from the choice of Chebychev samples.

For $t\in [-1.1,1.1]$, the relative interpolation errors are computed in the matrix Frobenius norm as
$
\frac{\|Q^*(t) - Q(t)\|_F}{\|Q(t)\|_F}
$, 
where $Q^*(t)$ denotes the manifold interpolant and $Q(t)$ is the reference solution.
The error curves are displayed in Fig. \ref{fig:Errors_analytic3_interp_n500p10_Cheby6_-1p1_1p1_invquad}.
%
%
%
%

%
%
The relative Frobenius errors are
\begin{center}
\begin{tabular}{l|l|l|l}
                     & Geo. interp. & RBF tan. interp. & Hermite interp.\\  \hline
Max. relative errors &        0.039 &    0.014     &  0.0007        \\ \hline
$L_2$ relative errors&        0.030 &    0.016     &  0.0005
\end{tabular}.
\end{center}
The (discrete) $L_2$-norm gives the integrated squared errors on a discrete representation of
the interval $[t_0, t_k]$ with a resolution of $100$ points.
%
%
%
%
\subsection{Hermite interpolation of a low-rank SVD}
\label{sec:CHSVDLowRank}
Next, we consider an academic example of a non-linear matrix function with fixed low rank.
The goal is to perform a quasi-cubic interpolation of the associated SVD.
As above, we construct a cubic matrix polynomial
\[ 
 Y(t) = Y_0 + t Y_1 + t^2Y_2 + t^3Y_3, \quad Y_i \in \R^{n\times r}, n=10,000, r=10
\]
with random matrices $Y_i$ with entries uniformly sampled from $[0,1]$ for $Y_0$
and from $[0,0.5]$ for $Y_1,Y_2, Y_3$.
Then, a second matrix polynomial is considered
\[ 
 Z(t) = Z_0 + t Z_1 + t^2Z_2, \quad Z_i \in \R^{r\times m}, r=10, n=300.
\]
Here, the entries of $Z_0$ are sampled uniformly from $[0,1]$
while the entries of $Z_1,Z_2$ are sampled uniformly from $[0,0.5]$.
The nonlinear low-rank matrix function is set as
\[
 W(t) = Y(t) Z(t) \in \R^{n\times m}.
\]
By construction, $W(t)$ is of fixed $\text{rank}(W(t)) \equiv r=10$ $\forall t$.
The low rank SVD
\[
 W(t) = U_r(t)\Sigma_r(t)V_r(t)^T, \quad U_r(t)\in St(n,r), V_r(t)\in St(m,r), \Sigma \in \R^{r\times r}
\]
is sampled at the two Chebychev nodes $t_0=0.0732, t_1=0.4268$ in the interval $[0.0, 0.5]$.\footnote{The Chebychev sampling is implemented as a standard in the program code that was written for the numerical experiments and is not of importance in this case.}
The Hermite sample data set
\[
 U_r(t_i), \dot U_r(t_i), \quad V_r(t_i), \dot V_r(t_i),\quad \Sigma_r(t_i), \dot \Sigma_r(t_i),\quad i=0,1,
\]
is computed with Alg. \ref{alg:SVDtruncdiff} of Appendix \ref{app:diffSVD}.
(To this end $V(t_i) = (V_r(t_i), V_{m-r}(t_i))\in O(m)$ is required.)

\begin{figure}[ht]
\centering
\includegraphics[width=0.6\textwidth]{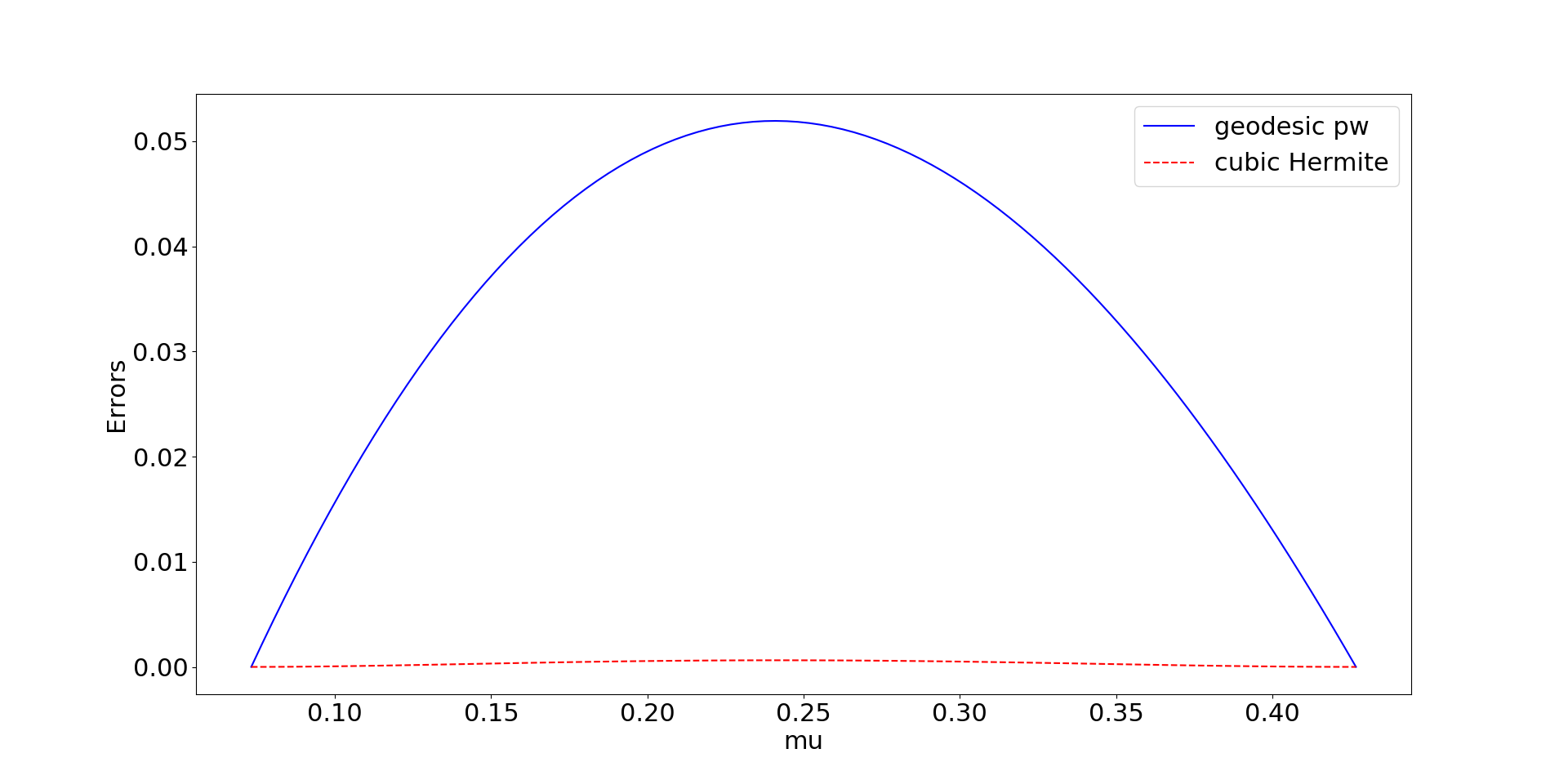}
\label{fig:Lowrank_interp_SVD_Terrors_n10000p10m300_Cheby2_-0_0p5}
\caption{Relative errors associated with the interpolation process of the low-rank SVD in terms of the Frobenius matrix norm.}
\end{figure}

\begin{remark}
Computing an analytic path of an SVD and thus a proper sample data set  is a challenge in its own right, see \cite{BunseGerstner1991}.
This is in part because of the inherent ambiguity of the SVD even in the case of mutually distinct singular values, where 
$W = U\Sigma V^T = (US) \Sigma (SV^T)$ 
for any orthogonal and diagonal matrix $S = \diag(\pm 1,\ldots, \pm 1)$, \cite[B.11, p. 334]{Higham:2008:FM}.
SVD algorithms from numerical linear algebra packages may return a different `sign-matrix' $S$ for the SVD of $W(t)$ and $W(s)$, even when 
$t$ and $s$ are close to each other. This introduces discontinuities in the sampled $U$ and $V$ matrices.
In the experiments performed in this work, we normalize the SVD as follows.
A reference SVD $U_0 \Sigma_0V_0^T = W(t_0)$ is computed. 
At each $t$, we compute an SVD $U_t \Sigma_tV_t^T$ and determine $S = \text{sign}(\diag(U_t^TU_0))$,
 where the sign-function is understood to be applied entry-wise on the diagonal elements.
Then, we replace $U_t \leftarrow U_tS$, $V_t \leftarrow V_tS$.
In the test cases considered here, this hands-on approach is sufficient to ensure a differentiable SVD computation.
In general, one has to allow for negative singular values to ensure differentiability, \cite{BunseGerstner1991}.
\end{remark}
For $t$ in the sampled range, the relative interpolation errors are computed in the Frobenius norm as
$
\frac{\|U^*(t)\Sigma^*(t)(V^*(t))^T - W(t)\|_F}{\|W(t)\|_F}
$, 
where $U^*(t)$, $\Sigma^*(t)$, $V^*(t)$ are the interpolants of the 
matrix factors of the low-rank SVD of $W(t)$ and $W(t)$ is the reference solution.
The relative errors are
\begin{center}
\begin{tabular}{l|l|l}
        & Geo. interp. & Hermite interp.\\  \hline
Max. relative errors &  0.0519 &  0.00063\\ \hline
$L_2$-norm of error data & 0.0225&  0.00024
\end{tabular}.
\end{center}
Fig. \ref{fig:Lowrank_interp_SVD_Terrors_n10000p10m300_Cheby2_-0_0p5} displays the error curves for the 
quasi-linear and the quasi-cubic Hermite interpolation approaches.

For the sake of completeness, we repeat this experiment but with selecting 
$p=U(\mu_i)$ as the center for the Riemannian normal coordinates.
Hence, the derivative data is mapped to $T_{U(\mu_i)}St(n,r)$ instead of $T_{U(\mu_{i+1})}St(n,r)$ 
and the tangent space interpolation curve is of the form
\begin{eqnarray*}
 \gamma(t) &=& a_0(t)\mathbf{0}_p + a_1(t) \Delta_q + b_0(t) \hat{v}_p + b_1(t) \hat{v}_q \subset T_pSt(n,r)
 \text{ instead of }\\
\gamma(t) &=& a_0(t) \Delta_p + a_1(t) \mathbf{0}_q + b_0(t) \hat{v}_p + b_1(t) \hat{v}_q \subset T_qSt(n,r),
\end{eqnarray*}
where $p =U(\mu_i) , q= U(\mu_{i+1})$.
This leads to virtually indistinguishable plots.
The maximum relative errors are 
$6.3023\cdot 10^{-4}$ ($q$-centered) vs. $6.3374\cdot 10^{-4}$ ($p$-centered).
The $L_2$-norms of the relative errors are
$2.3819\cdot 10^{-4}$ ($q$-centered) vs. $2.3954\cdot 10^{-4}$ ($p$-centered).

\begin{figure}[ht]
\centering
\includegraphics[width=0.6\textwidth]{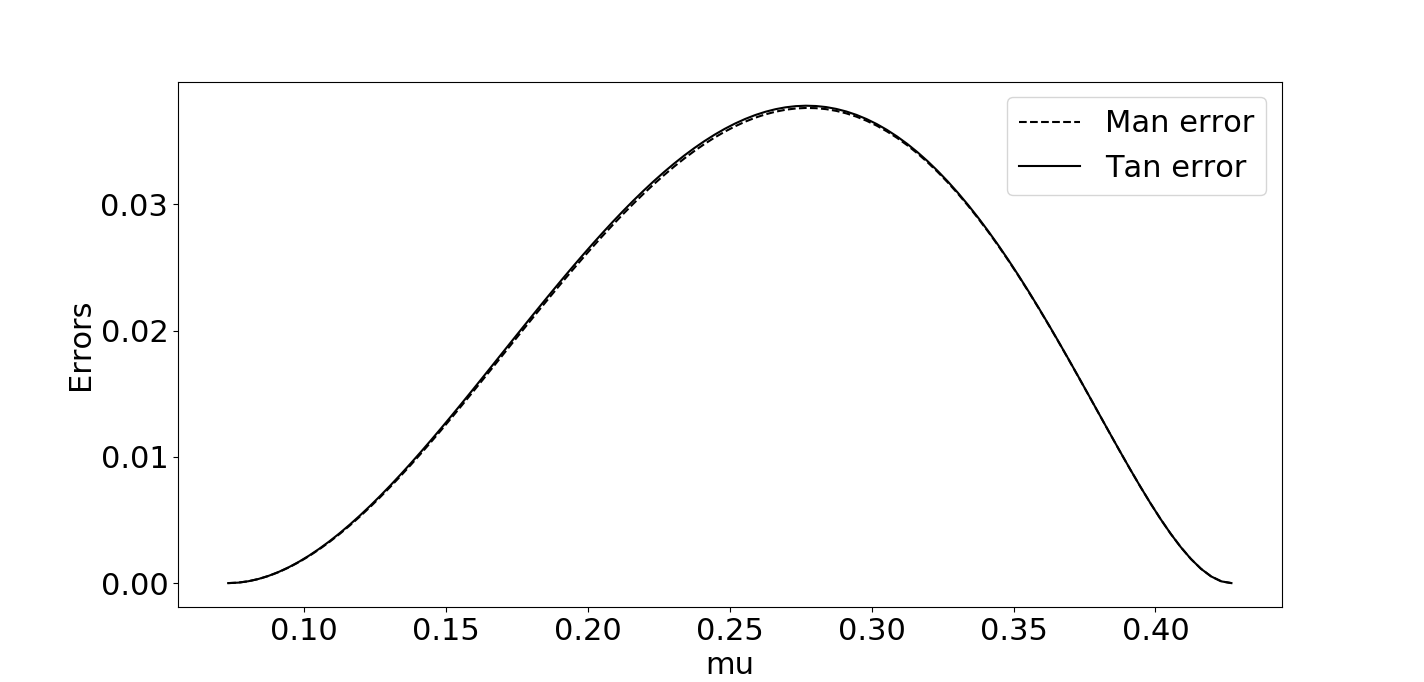}
\label{fig:Lowrank_interp_SVD_ManTan_errors_n10000p10m300_Cheby2_-0_0p5}
\caption{Interpolation of the $U$-factor of the SVD data in Section \ref{sec:CHSVDLowRank}.
Absolute Hermite interpolation errors in terms of the Riemannian metric on the tangent space (Tan error) and as measured by the 
Riemannian distance function on the manifold (Man error).}
\end{figure}

Recall that the local cubic Hermite interpolation scheme works in essence by performing 
Hermite interpolation in a selected tangent space and subsequently mapping the result to the manifold.
For the $U$-factor interpolation featured in the above example,
Fig. \ref{fig:Lowrank_interp_SVD_ManTan_errors_n10000p10m300_Cheby2_-0_0p5} shows the absolute
interpolation errors of the tangent space data in the {\em canonical Riemannian metric}
together with interpolation errors of final manifold data in terms of the {\em Riemannian distance}.
The manifold errors are very close to the tangent space errors but are actually slightly smaller,
inspite of the additional downstream translation of the tangent space interpolants 
to the manifold via the Riemannian exponential, which is an additional source of numerical errors.
This is in line with Theorem \ref{thm:ErrorMain}, since the Stiefel manifold features positive sectional curvature. 
%
%
%
\subsection{Hermite interpolation of the left singular values of non-linear function snapshots}
\label{sec:CHSVD_fun}
In the next experiment, we consider the SVD of discrete snapshots of
a nonlinear multi-parameter function.
To this end, define 
\begin{eqnarray*}
f&:&[0,1]\times [0,2]\times [1,4] \rightarrow \R, (x,t,\mu)\mapsto x^t \sin(\frac{\pi}{2}\mu x) \mbox{ and}\\
F&:&[0,1]\times [0,2]\times [1,4]\rightarrow \R, 
	(x,t,\mu)\mapsto \frac{f(x,t,\mu)}{\|f(\cdot, t,\mu)\|_{L_2}},
\end{eqnarray*}
where 
$
  \langle f_1, f_2 \rangle_{L_2} = \int_0^{1}{f_1(x)f_2(x)dx}
$
and $\|\cdot\|_{L_2} = \sqrt{\langle \cdot, \cdot \rangle_{L_2}}$ on $L_2([0,1])$.
We will discretize $F$ in $x$, take function `snapshots' at selected values of $t$
and eventually Hermite interpolate the left singular vectors of the discrete snapshot matrices with respect to $\mu$.
The partial derivative of $F$ by $\mu$ is
\begin{equation}
	\partial_{\mu}F(x,t,\mu)  =
	\frac{1}{\|f(\cdot, t,\mu)\|_{L_2}} \partial_{\mu}f(x,t,\mu)
	- \frac{\langle f(\cdot, t,\mu), \partial_{\mu} f(\cdot, t,\mu)\rangle_{L_2}}
	{\|f(\cdot, t,\mu)\|_{L_2}^3}f(x,t,\mu).
\end{equation}
For the spatial discretization, we use an equidistant decomposition of the unit interval, $0=x_1, x_2,\ldots, x_n=1$, $n=1001$.
Then, we take $r=6$ function snapshots in $t$ at 
time instants $t= 1.0,\hspace{0.1cm}   1.6,\hspace{0.1cm}  2.2,\hspace{0.1cm}    2.8, \hspace{0.1cm}   3.4, \hspace{0.1cm}   4.0$.
In this way, a $\mu$-dependent snapshot matrix function 
$Y(\mu) := (F(x,t_1,\mu),\ldots,F(x,t_6,\mu))\in  \R^{n\times r} = \R^{1001\times 6}$
with SVD
\[
 Y(\mu) = U(\mu)\Sigma(\mu) V(\mu)^T, \quad U(\mu)\in St(n,r), \Sigma(\mu)\in\R^{r\times r}, V(\mu)\in O(r)
\]
is obtained.
\begin{figure}[ht]
\centering
\includegraphics[width=1.0\textwidth]{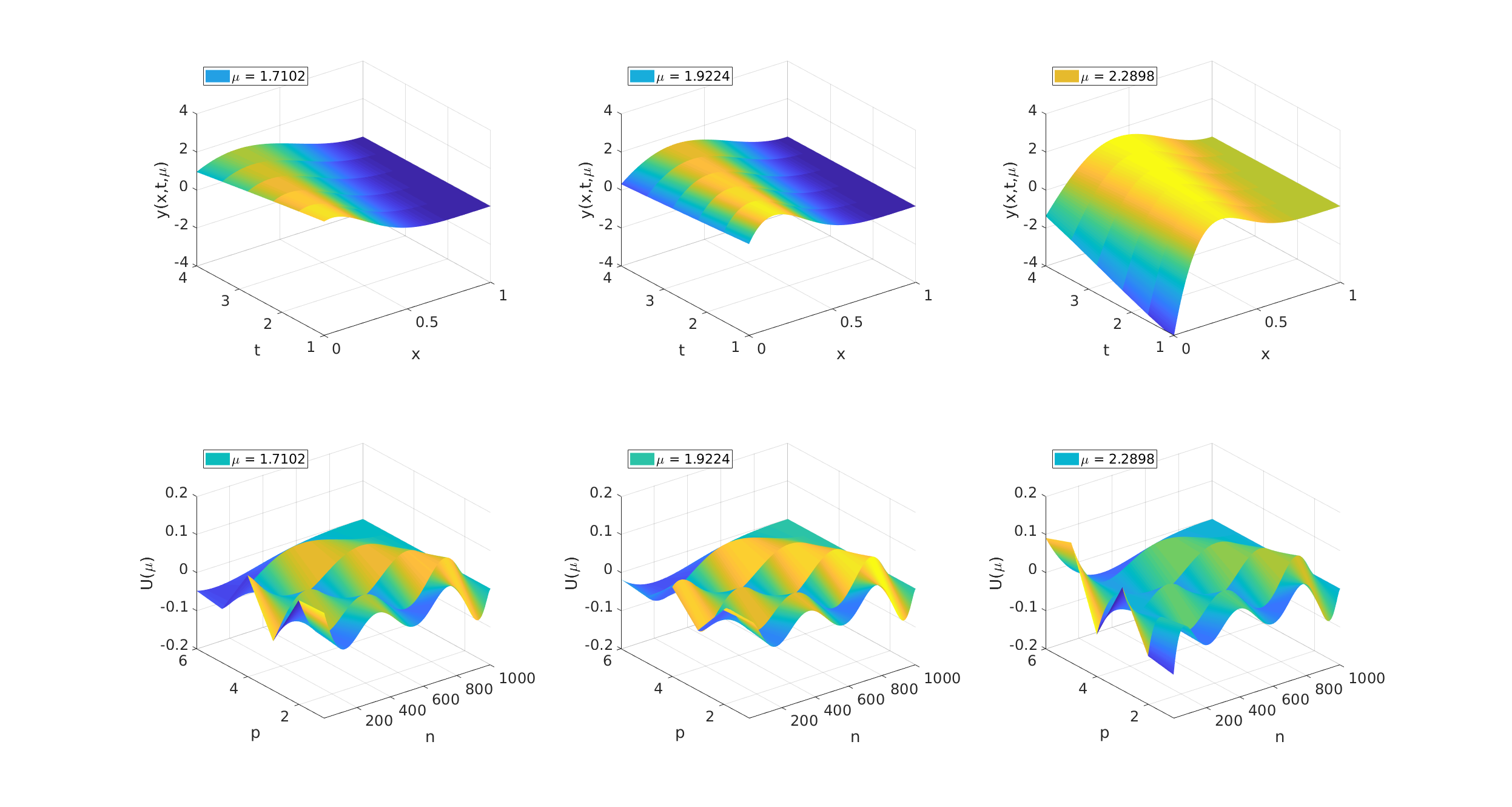}
\caption{Sample data featured in Section \ref{sec:CHSVD_fun}. Upper row: snapshot matrices $Y(\mu)$, lower row: corresponding
left singular vector matrices $U(\mu)$.
}
\label{fig:sample_set_analytic2}
\end{figure}
%
Fig. \ref{fig:sample_set_analytic2} displays the snapshot matrices at some selected parameter values,
together with the associated left singular value matrices.
For $\mu\in [0,2]$, the values of $f$ are non-negative and so are all entries in the corresponding snapshot matrices.
Beyond $\mu = 2.0$, negative entries arise in the snapshot vectors.
Fig. \ref{fig:analytic2_interp_n1001p6_6cheby_1p7_2p3_Sigma_p_Plot} tracks the smallest singular value $\sigma_r(\mu)$ of
the snapshot matrices $Y(\mu)$ for $\mu \in [1.7, 2.3]$.
A substantial non-linear change in $\sigma_r(\mu)$ is apparent around the value of $\mu=2.0$.
This makes SVD interpolation beyond the parameter location $\mu=2.0$ a challenging problem.

We sample the left singular value matrix $U(\mu)$ together with the derivative $\dot{U}(\mu)$ 
at 6 Chebychev samples in the interval $[1.7,2.3]$.
%
\begin{figure}[ht]
\centering
\includegraphics[width=0.6\textwidth]{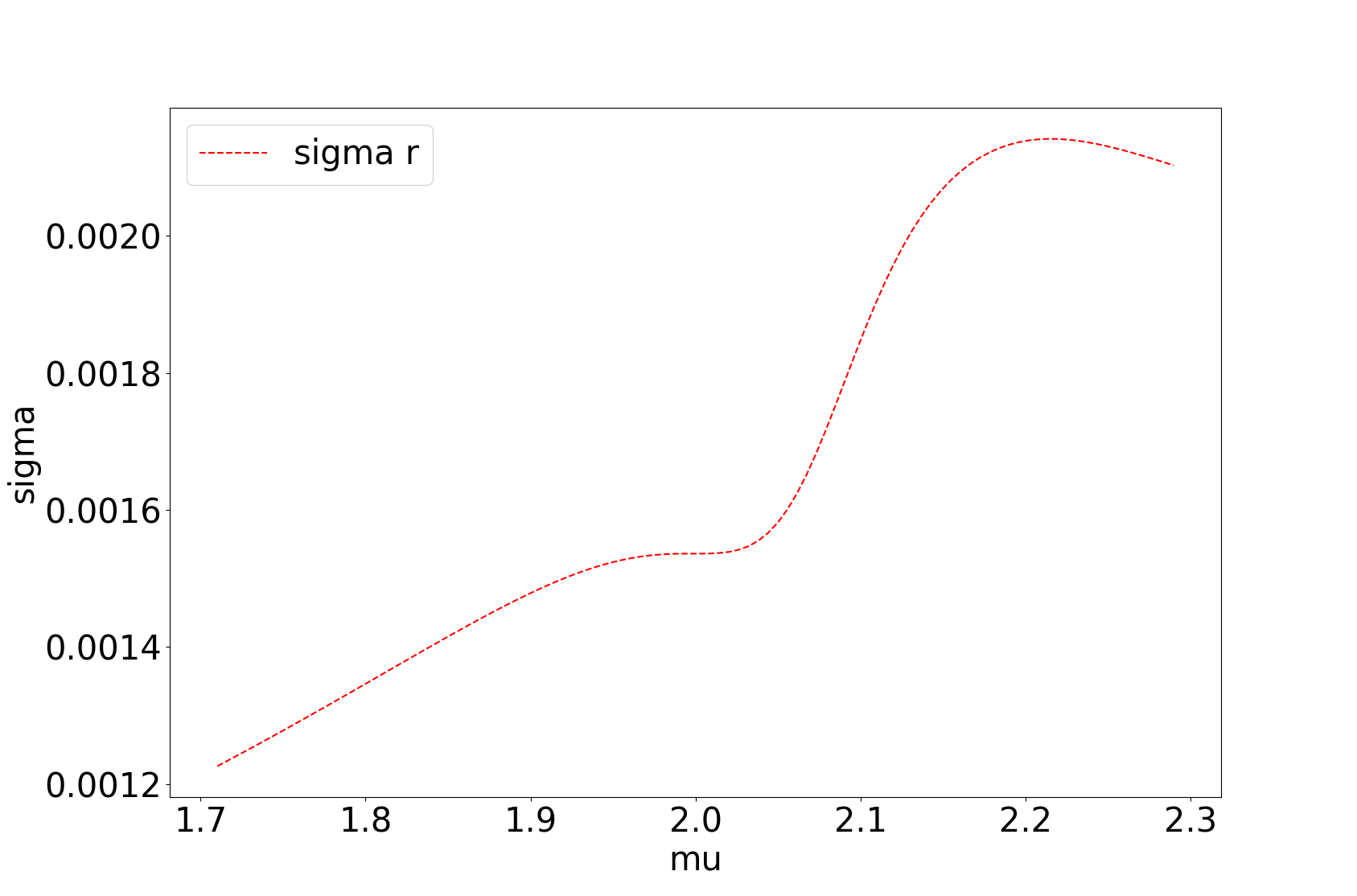}
\caption{Corresponding to Section \ref{sec:CHSVD_fun}. The smallest singular value $\sigma_r(\mu)$ of the $\mu$-dependent SVD $Y(\mu) = U(\mu)\Sigma(\mu) V(\mu)^T$
in the range $\mu\in [1.7, 2.3]$.
}
\label{fig:analytic2_interp_n1001p6_6cheby_1p7_2p3_Sigma_p_Plot}
\end{figure}
%
As in Section \ref{sec:CHQR}, we juxtapose the results of quasi-linear interpolation,
tangent space interpolation and quasi-cubic Hermite interpolation.
For $\mu$ in the sampled range, the relative interpolation errors are computed in the Frobenius norm as
$
\frac{\|U^*(\mu) - U(\mu)\|_F}{\|U(\mu)\|_F}
$, 
where $U^*(\mu)$ denotes the manifold interpolant and $U(\mu)$ is the reference solution.
The error curves are displayed in Fig. \ref{fig:analytic2_errors}.
According to the figure, the tangent space interpolation method fails to interpolate
the samples at the first two parameter locations $\mu_0 \approx 1.7102,\mu_1\approx 1.7879$.
This is explained as follows. In the tangent space interpolation method,
all the Stiefel samples $U(\mu_i)$ are mapped to the tangent space
attached at $U(\mu_3)$, $\mu_3 \approx 2.0776$ via
$\Delta(\mu_i) = \Log_{U(\mu_3)}^{St}(U(\mu_i)$.\footnote{Note that the base point $\mu_3$ happens to lie beyond the `$2.0$-threshold value', after which negative function values appear,
while $\mu_0, \mu_1< 2.0$.}
It turns out that the Riemannian Stiefel logarithm is not well-defined for $i=0,1$.
Put in different words, $U(\mu_0)$ and $U(\mu_1)$ are too far from $U(\mu_3)$ to be mapped to $T_{U(\mu_3)}St(n,r)$
by the the Stiefel $\log$-algorithm.
%
%
\begin{figure}[ht]
\centering
\includegraphics[width=0.8\textwidth]{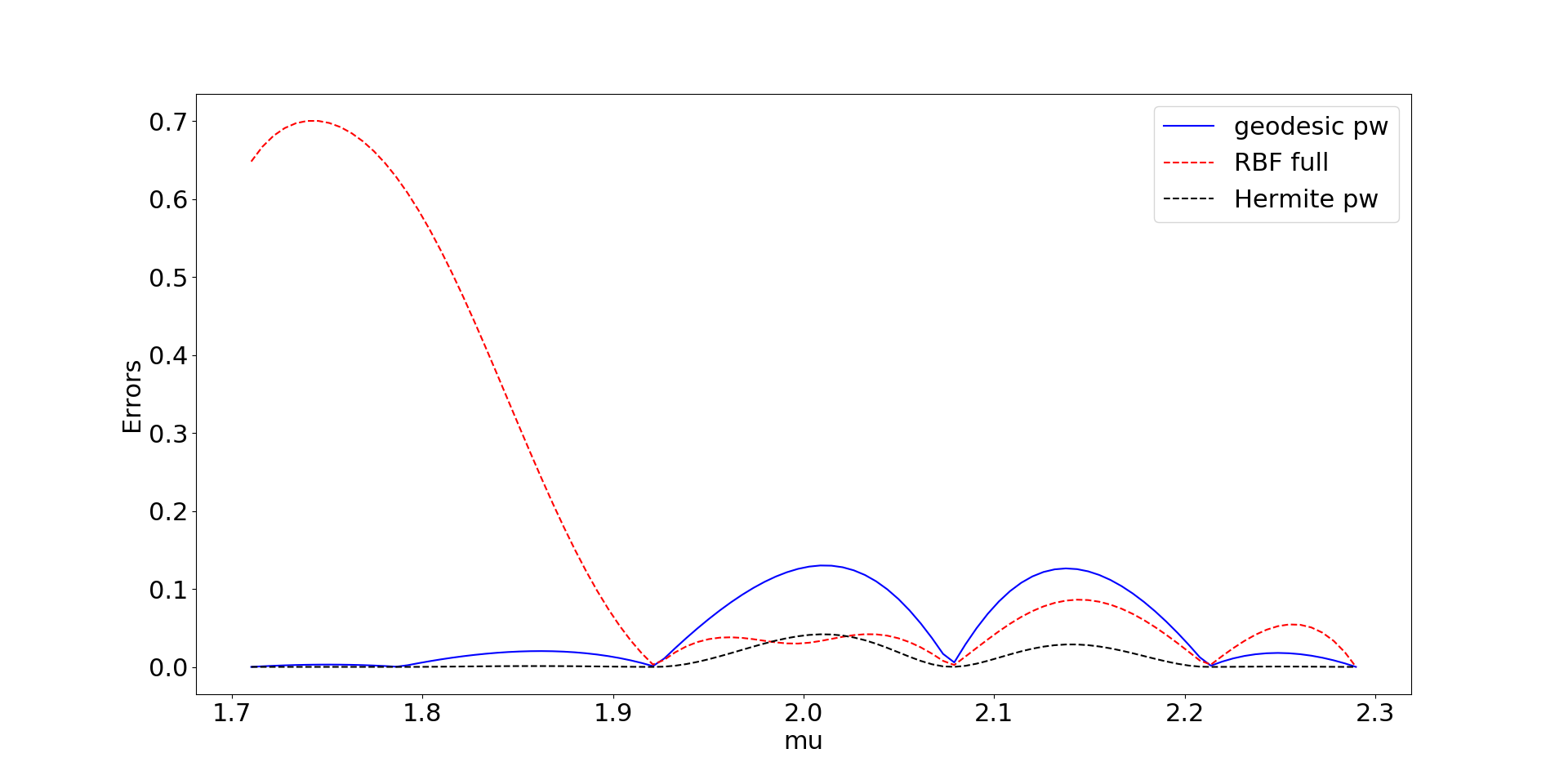}
\caption{Sample data
}
\label{fig:analytic2_errors}
\end{figure}
%

%
The relative Frobenius errors are
\begin{center}
\begin{tabular}{l|l|l|l}
                     & Geo. interp. & RBF tan. interp. & Hermite interp.\\  \hline
Max. relative errors &  0.1301      &    0.7003        &  0.0418        \\ \hline
$L_2$-norm of error data&  0.0501      &    0.2336        &  0.0123
\end{tabular}.
\end{center}
%
%
%
%
%
%
%
\subsection{Parametric dimension reduction}
\label{sec:BS_interp}
SVD interpolation may be used for parametric dimension reduction.
Consider samples $Y_0,\ldots Y_k$ of a matrix curve $\mu\mapsto Y(\mu)\subset \R^{n\times m}$.
Suppose that instead of the original data, only a low-rank SVD approximation 
is stored $Y_i \approx  U_i\Sigma_iV_i^T$, where $U_i\in \R^{n\times r}, \Sigma_i\in \R^{r\times r}, V_i\in\R^{m\times r}$
for a fixed $r\ll \min\{n,m\}$.
Then a low-rank approximant to any $Y(\mu)$ in the sampled range can be obtained via interpolating the SVD data.

In this section, we apply this approach to an application from computational option pricing.
The value function $y(T,S;K,r,\sigma)$ that gives the fair price for a European call option is determined via the {\em Black-Scholes-equation} \cite{black73tpo},
\begin{align*}
y_t(t,S)&=\frac12\sigma^2S^2y_{SS}(t,S)+rSy_S(t,S)-ry(t,S),\quad S\geq0,\quad 0<t\leq T,\\
y(T,S)&=\max\{S-K,0\},\quad S\geq0.
\end{align*}
This is a parabolic PDE that depends on time $t$, the stock value $S$, and a number of additional system parameters, namely the strike price $K$, the interest rate $r$, the volatility $\sigma$ and the exercise time $T$.
In this experiment, we consider a fixed interest rate of $r=0.01$ and an exercise time of $T=2$ units.
The dependency on the underlying $S\in [50,150]$ is resolved via a discretization of the interval by equidistant steps of $\Delta S=0.01$, while the strike price $K\in[30,170]$ is discretized in steps of $\Delta K=1$.
Eventually, the volatility $\sigma$ will act as the interpolation parameter.
Hermite interpolation requires the option price $y$ as well as its derivative $\partial_{\sigma}y$,
in economics referred to as the `vega' of set of the `greeks'. A similar test case was considered in \cite{Zimmermann2018}.

\begin{figure}[ht]
\centering
\includegraphics[width=0.8\textwidth]{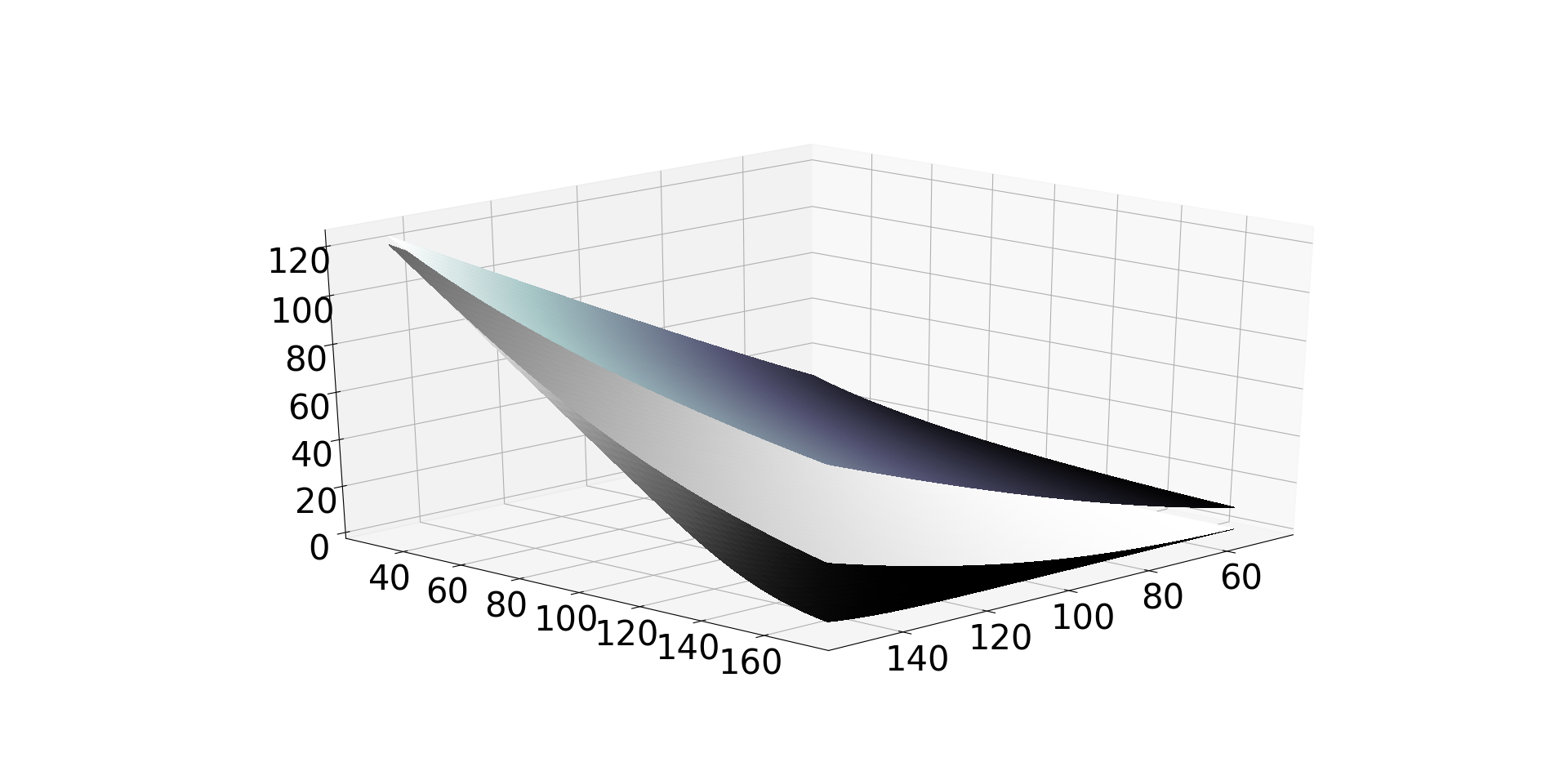}
\label{fig:rom_OPTION}
\caption{Sampled data Black-Scholes solution data sets $Y(0.1), Y(0.4), Y(0.9)$}
\end{figure}

\begin{figure}[ht]
\centering
\includegraphics[width=0.9\textwidth]{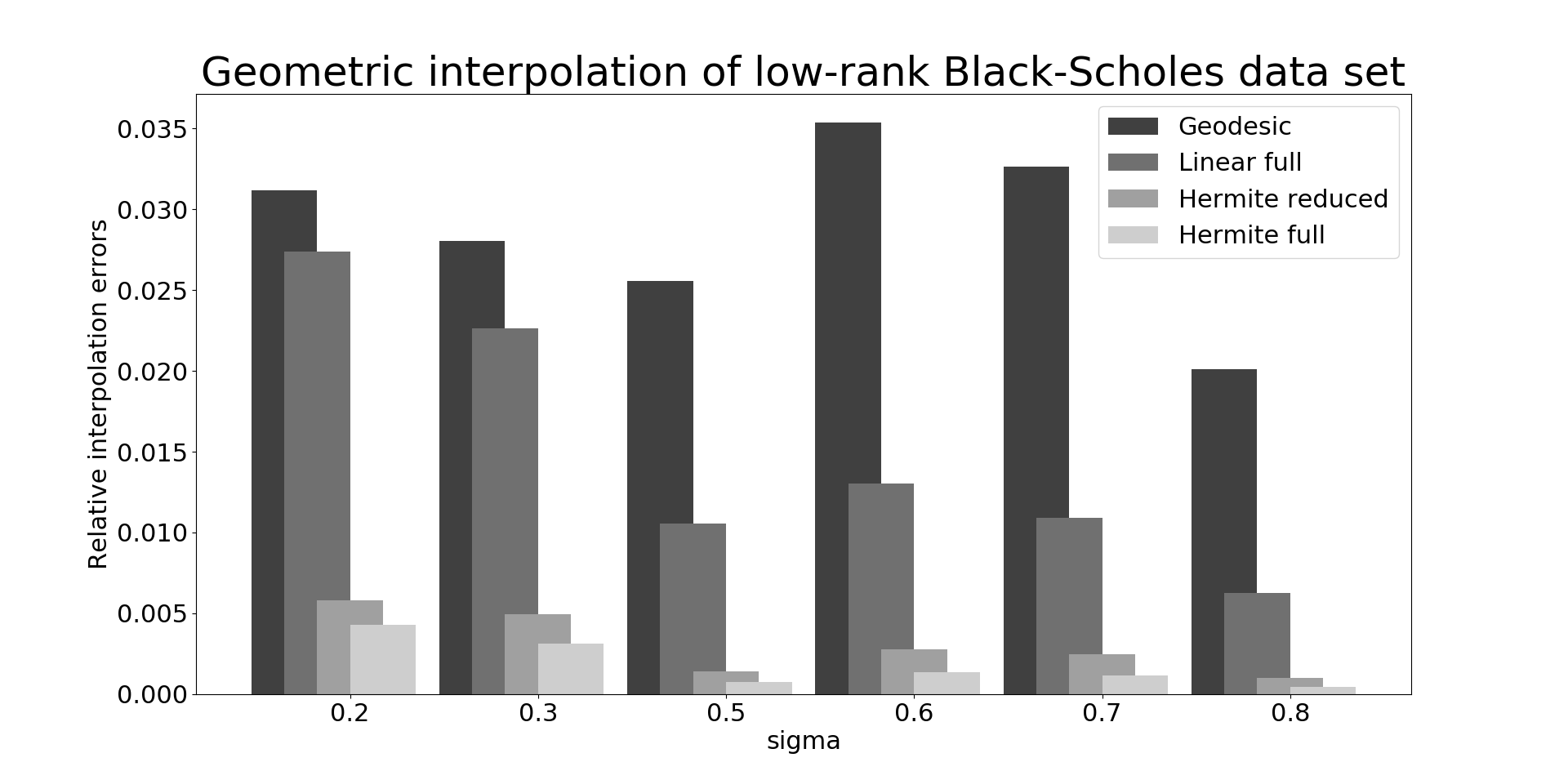}
\label{fig:Hermite_OPTION}
\caption{Relative interpolation errors associated with the Black-Scholes test case. Compressed, interpolated SVD obtained via geodesic interpolation
 and Hermite interpolation, respectively, vs. uncompressed standard linear and Hermite interpolation.}
\end{figure}

The Black-Scholes equation for a single underlying has a closed-form solution. Yet, here, we will approach it via a numerical scheme in order to mimic the corresponding procedure for real-life problems.
Application of a finite volume scheme to the Black-Scholes PDE yields snapshot matrices
$$Y(\sigma)=\left(Y_{s,k}(\sigma)\right)_{\substack{s=1,\dots,10001\\k=1,\dots,141}},\quad \partial_\sigma Y(\sigma)_{\substack{s=1,\dots,10001\\k=1,\dots,141}},$$
  for $\sigma\in[0.1,0.2,\dots,1.0]$.
The computation time for each data pair $Y(\sigma), \partial_\sigma Y(\sigma)$ is ca. 17min on a standard laptop computer.
%
For each sampled snapshot matrix $Y(\sigma)\in \R^{n\times m}$, $n= 10001$, $m = 141$,
an SVD is performed and is truncated to the $r=5$ dominant singular values/singular vector triples.
This yields a compressed representation
$Y(\sigma) = U(\sigma)\Sigma(\sigma) V^T(\sigma)$,
with  $U(\sigma) \in St(n,r)$, $\Sigma(\sigma) \in diag(r,r)$,
$V(\sigma) \in St(m,r)$ and consumes ca. $0.1s$ on a laptop computer.
The {\em relative information content} is ric$(r) = \frac{\sum_{j=1}^r \sigma_j^2}{\sum_{k=1}^m \sigma_k^2} \geq 0.99999$.
The storage requirements for a $(10001\times 141)$-matrix are 11.3MB,
all the low-rank SVD factors truncated to $r=5$ require a total of 0.4MB of disk space, which is ca. $3.5\%$ of the uncompressed representation.
We sample full solution data sets
$Y(\sigma)$, $\partial_\sigma Y(\sigma)$ at $\sigma \in \{0.1, 0.4, 0.9\}$.
The sample data sets are displayed in Fig. \ref{fig:rom_OPTION}.
%
%
%
In order to assess the approximation accuracy, we interpolate at
 $\sigma^*\in\{0.2,0.3,\hspace{0.1cm}0.5, 0.6,0.7,0.8\}$
and compute the relative Frobenius norm errors of the interpolated
low-rank SVD representation $U(\sigma^*)\Sigma(\sigma^*)V(\sigma^*)^T$
with respect to the exact full rank data matrix $Y(\sigma^*)$.
We compare the quasi-linear geodesic interpolation (w/o derivative data)
to the cubic Hermite approach.
The errors are displayed in the bar plot Fig. \ref{fig:Hermite_OPTION}.
For completeness, we also include the results of standard linear and Hermite
interpolation on the data set of the full, uncompressed matrices $Y_i$, where no special geometric structure needs to be addressed.
%
The values underlying the bar plot are
\begin{footnotesize}
\begin{verbatim}
sigma         0.2      0.3      0.5      0.6      0.7      0.8
Geo.         [0.031178 0.028025 0.025577 0.035379 0.032636 0.0201280]
Linear full  [0.027382 0.022640 0.010562 0.013015 0.010908 0.0062528]
Hermite      [0.005791 0.004969 0.001416 0.002799 0.002482 0.0009872]
Hermite full [0.004298 0.003142 0.000734 0.001374 0.001156 0.0004374]
\end{verbatim}
\end{footnotesize}

Mind that the errors for the geodesic and the Hermite low rank interpolation
include both the effects of interpolation and data reduction.
Even though the data sets underwent a substantial reduction in dimension,
the relative errors are of a comparable order of magnitude.
For better judging the results, we note that the relative error between
the two consecutive samples $Y(0.1)$, $Y(0.4)$ and $Y(0.4)$, $Y(0.9)$
are $27\%$ and $39\%$, respectively.

The technique could be used as a reduced online storage scheme:
store the truncated (Hermite) SVD data at some selected sample locations;
interpolate, when a prediction at any in-between location is required online.
%
%
%
%
%
%
%
%
%
%
%
%
%
\section{Conclusions and final remarks}
\label{sec:conclusions}
We have presented an elementary, general approach to Hermite interpolation on Riemannian manifolds
that is applicable to practical problems, whenever algorithms to compute the Riemannian exp and log mappings are available.
While our focus was on the manifold counterpart of local cubic Hermite interpolation, 
the method is flexible and may be combined with any Hermite method that is linear in the sample data.
In fact, only the coefficient functions $a_0, b_0, b_1$ in \eqref{eq:cubicHermiteMnf} need to be replaced,
no additional changes are necessary.
Moreover, combinations of Hermite and Lagrange methods are straightforward generalizations.

In addition, we have exposed a relation between the sectional curvature of the manifold in question the data processing errors, that arise for computations in Riemannian normal coordinates.

As an example, Hermite interpolation of Stiefel data was discussed in more detail.
From the observations in the numerical experiments, the main practical constraint on the sampled data is that two consecutive samples be close enough 
so that the Riemannian Stiefel logarithm is well-defined.
As a rule of thumb, if the data points are close enough 
so that the Riemannian log algorithm converges, then
the Hermite interpolation method provides already quite accurate results.

The method constructs piece-wise cubic manifold splines between data points $p_{i}$ and $p_{i+1}$
in terms of normal coordinates centered at $p_{i+1}$. Thus, it is not symmetric in the sense that 
computations in normal coordinates centered at $p_i$ might lead to different results.
Yet, in the numerical experiments, these effects prove to be negligible.
\section*{Acknowledgments}
The data set featured in Section \ref{sec:BS_interp} was kindly provided by my colleague Kristian Debrabant from 
the Department for Mathematics and Computer Science (IMADA), SDU Odense.
\appendix
%

%
%
%
\section{The basic cubic Hermite coefficient polynomials}
\label{app:classicHermite}
\begin{figure}[ht]
\centering
\includegraphics[width=0.8\textwidth]{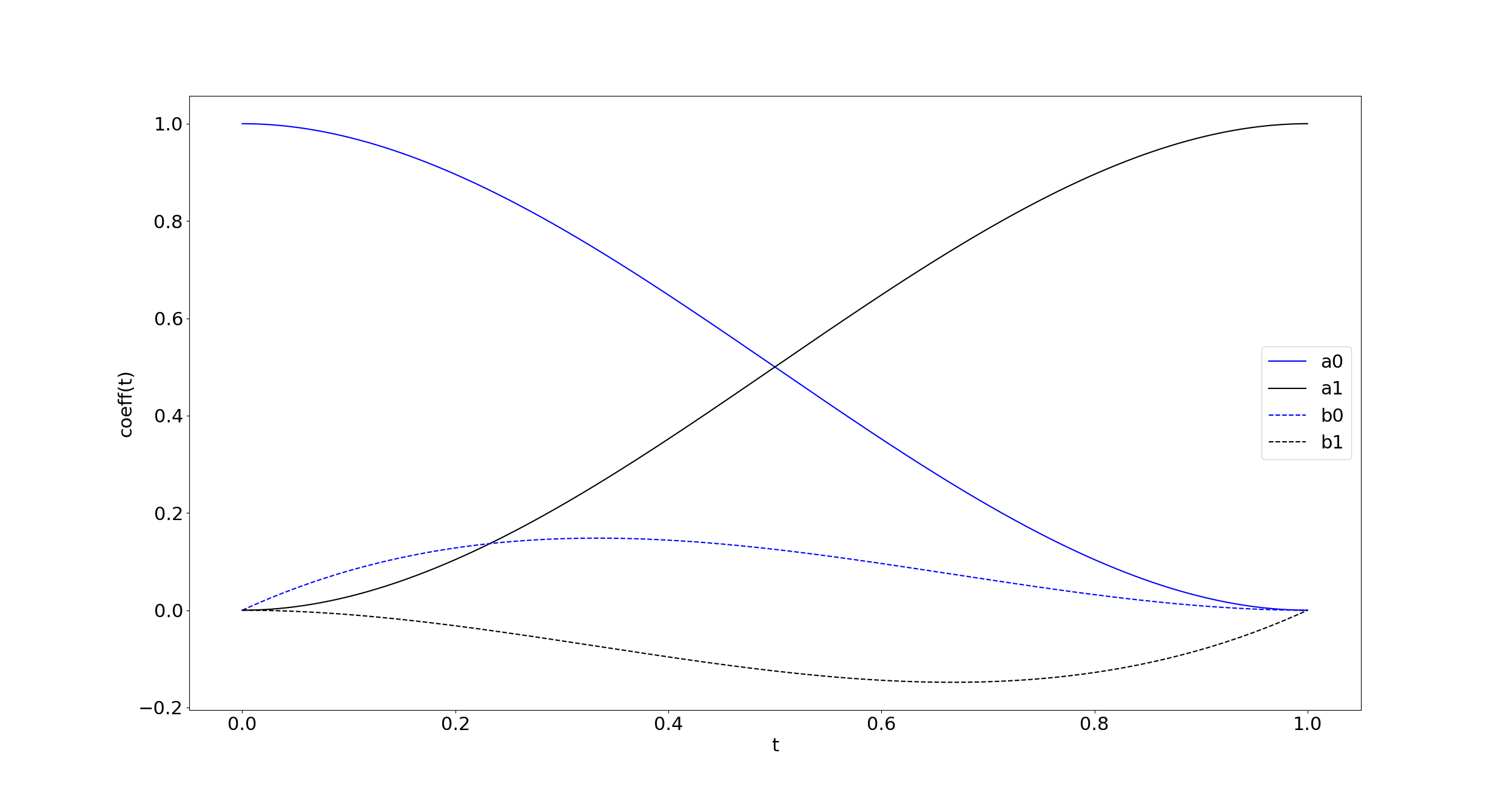}
\caption{The classical cubic Hermite coefficient functions for $t_0=0, t_1=1$,
$a_0(t) = h_{1000}(t)= 1-3t^2 + 2t^3$,
$a_1(t) = h_{0010}(t)= 3t^2 - 2t^3$, 
$b_0(t) = h_{0100}(t)=t-2t^2+t^3$, 
$b_1(t) = h_{0001}(t)=t^3-t^2$.
}
\label{fig:Cubic_Hermite_coeffs}
\end{figure}

The coefficient functions in \eqref{eq:VectorCubicHermite}
are the cubic Hermite polynomials that are uniquely defined by
\[
 \begin{array}{c|c|c|c|c|c}
f & f(t_0)&f'(t_0) &f(t_1) & f'(t_1) \\
       \hline
a_0 & 1   & 0      & 0     & 0   \\
a_1 & 0   & 0      & 1     & 0   \\
b_0 & 0   & 1      & 0     & 0   \\
b_1 & 0   & 0      & 0     & 1   \\
 \end{array}.
\]
The explicit cubic coefficient functions are
\begin{eqnarray}
\label{eq:a0}
 a_0(t) &=& 1 - \frac{1}{(t_1-t_0)^2}(t-t_0)^2 +\frac{2}{(t_1-t_0)^3}(t-t_0)^2(t-t_1),\\
 \label{eq:a1}
 a_1(t) &=& \frac{1}{(t_1-t_0)^2}(t-t_0)^2 -\frac{2}{(t_1-t_0)^3}(t-t_0)^2(t-t_1),\\
 \label{eq:b0}
 b_0(t) &=&  (t-t_0) - \frac{1}{(t_1-t_0)}(t-t_0)^2 +\frac{1}{(t_1-t_0)}(t-t_0)^2(t-t_1),\\
\label{eq:b1}
 b_1(t) &=&  \frac{1}{(t_1-t_0)^2}(t-t_0)^2(t-t_1),
\end{eqnarray}
and are displayed in  Fig. \ref{fig:Cubic_Hermite_coeffs} for $t_0=0, t_1 =1$.
Since on manifolds, we work exclusively in the setting, where $q=0$,
the coefficient $a_1(t)$ drops out in \eqref{eq:VectorCubicHermite}.

%
%

\begin{figure}[ht]
\centering
\includegraphics[width=0.7\textwidth]{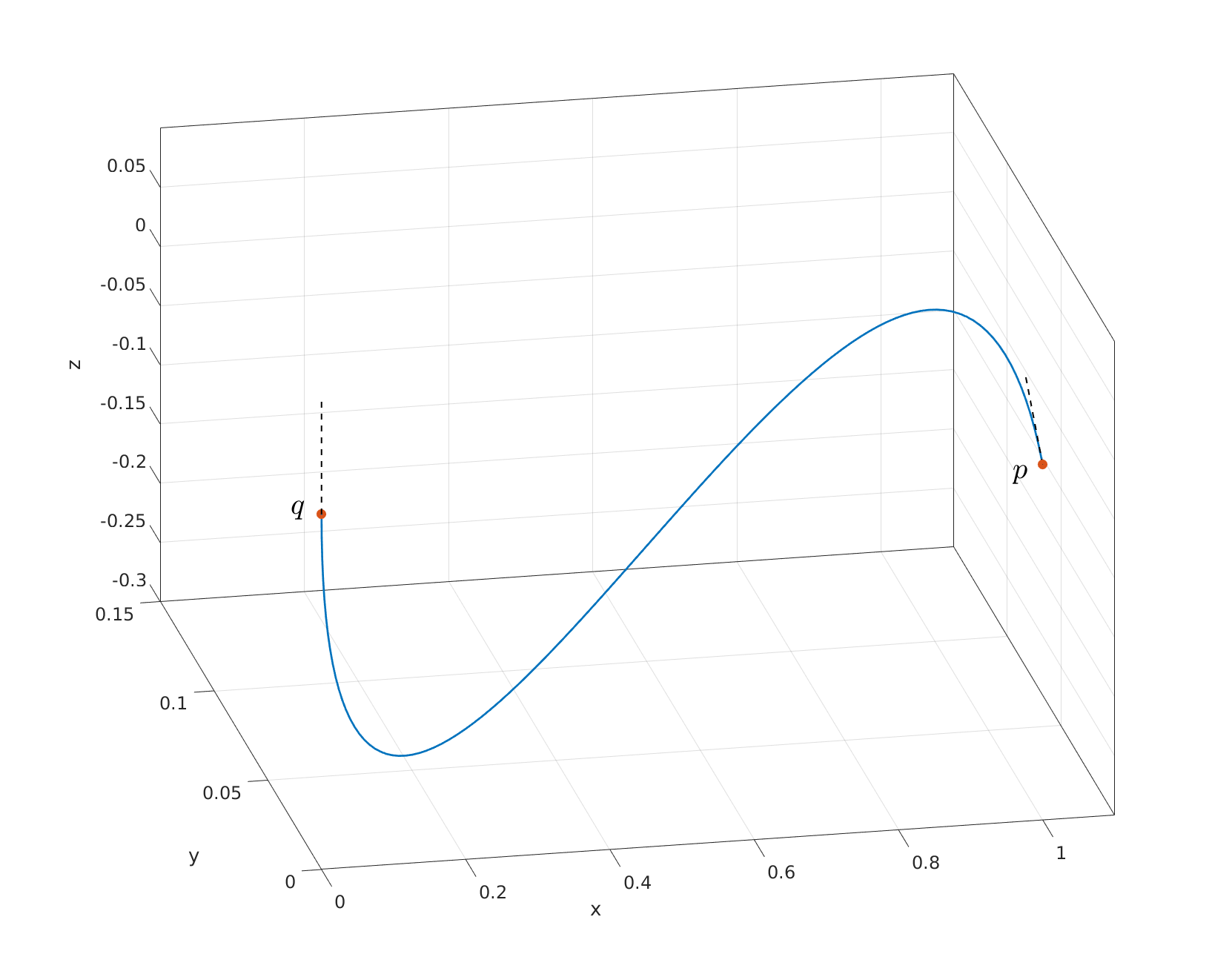}
\caption{Cubic Hermite spline in $\R^3$ starting from $p=(1,0,0)$ with velocity $v_p=(0.5,0.5,0)$
and ending in $q = (0,0,0)$ with velocity $v_q=(0,0,1)$. (Velocity directions indicated by the dashed lines.)
}
\label{fig:Cubic_Hermite_spline}
\end{figure}

%
%
%
Fig. \ref{fig:Cubic_Hermite_spline} shows the spatial cubic Hermite spline \eqref{eq:VectorCubicHermite} that connects the points 
$p=(1,0,0), q = (0,0,0)\in\R^3$ with a prescribed start and terminal velocity of
$v_p=(0.5,0.5,0)$ and $v_q=(0,0,1)$, respectively.
%
%
%
%
%
%
\section{Differentiating the QR-decomposition}
\label{app:diffQR}
Let $t\mapsto T(t)\in \R^{n\times r}$ be a differentiable matrix function with Taylor expansion
$T(t_0+ h) = T(t_0) + h \dot{T}(t)$.
Following \cite[Proposition 2.2]{WalterLehmannLamour2012}, the QR-decomposition is characterized via the following set of matrix equations.
\[
 T(t)  = Q(t)R(t),\quad Q^T(t)Q(t) = I_r, \quad 0 = P_L\odot R(t).
\]
In the latter, 
$P_L = 
                    \left(
                     \begin{smallmatrix}
                       0      & 0      &\cdots & 0\\
                       1      & \ddots &\ddots & 0\\
                       \vdots & \ddots &\ddots &\vdots \\
                       1      & \cdots &1      & 0
                      \end{smallmatrix}
                    \right)
$
and `$\odot$' is the element-wise matrix product so that
$P_L\odot R$ selects the lower triangle of the square matrix $R$.
For brevity, we write $T= T(t_0), \dot T = \frac{d}{dt}\big\vert_{t=t_0}T(t)$, likewise for $Q(t)$, $R(t)$.
By the product rule
\[
 \dot{T} = \dot{Q} R + Q\dot{R}, \quad 0 =\dot{Q}^TQ + Q^T\dot{Q}, \quad 0  = P_L\odot \dot{R}.
\]
According to \cite[Proposition 2.2]{WalterLehmannLamour2012}, the derivatives $\dot{Q}, \dot{R}$
can be obtained from Alg. \ref{alg:QRdiff}.
The trick is to compute $X = Q^T\dot{Q}$ first and then use this to compute
$\dot{Q} = QQ^T\dot{Q} + (I-QQ^T)\dot{Q}$ by exploiting that $Q^T\dot{Q}$ is skew-symmetric
and that $\dot{R}R^{-1}$ is upper triangular.
\begin{algorithm}
\caption{Differentiating the QR-decomposition, \cite[Proposition 2.2]{WalterLehmannLamour2012}}
\label{alg:QRdiff}
\begin{algorithmic}[1]
  \REQUIRE{matrices $T, \dot{T}\in\R^{n\times r}$, (compact) QR-decomposition $T=QR$.}
  \STATE{ $L:= P_L\odot(Q^T\dot{T}R^{-1})$}
  \STATE{ $X = L-L^T$} \hfill \COMMENT{Now, $X = Q^T\dot{Q}$}
  \STATE{ $\dot{R} =  Q^T\dot{T} - XR$}
  \STATE{ $\dot{Q} = (I-QQ^T)\dot{T}R^{-1} + QX$}
  \ENSURE{$\dot{Q}, \dot{R}$}
\end{algorithmic}
\end{algorithm}
\section{Differentiating the singular value decomposition}
\label{app:diffSVD}
Let $m\leq n\in \N$ and suppose that $t \mapsto Y(t)\in \R^{n\times m}$ is a differentiable matrix curve around $t_0\in \R$.
If the singular values of $Y(t_0)$ are mutually distinct, then
the singular values and both the left and the right singular vectors 
depend differentiable on $t \in [t_0 -\delta t, t_0+\delta  t]$
for $\delta t$ small enough. 
This is because the associated symmetric eigenvalue problem $Y^T(t)Y(t) = V(t)\Lambda(t) V^T(t)$
is differentiable under these (and more relaxed) conditions, \cite{Alekseevsky1998}.

Let $ t \mapsto Y(t) = U( t)\Sigma( t) V( t)^T \in \R^{n\times m}$, 
where $U(t)\in St(n,m)$, $V( t)\in O(m)=St(m,m)$ and $\Sigma( t)\in\R^{m\times m}$ diagonal and positive definite.
Let $u_j$ and $v_j$, $j=1,\ldots,m$ denote the columns of $U( t_0)$ and $V( t_0)$, respectively.
For brevity, write $Y = Y(t_0), \dot{Y} = \frac{d}{dt}\big\vert_{t=t_0}Y(t)$, likewise for the other matrices that feature in the SVD.

\begin{algorithm}
\caption{Differentiating the SVD}
\label{alg:SVDdiff}
\begin{algorithmic}[1]
  \REQUIRE{matrices $Y, \dot{Y}\in\R^{n\times m}$, (compact) SVD $Y = U\Sigma V^T$.}
  \STATE{ $\dot{\sigma}_j = (u_j)^T \dot{Y} v_j\mbox{ for }j = 1,\ldots, m$}
  \STATE{ $\dot{V} =V \Gamma, \mbox{ where } \Gamma_{ij} =
            \left\{
                \begin{array}{ll}
                 \frac{ \sigma_i (u_i^T \dot{Y} v_j) +  \sigma_j(u_j^T\dot{Y} v_i) }{(\sigma_j + \sigma_i)(\sigma_j - \sigma_i)}, & i\neq j\\
                0,          & i=j
                \end{array}
            \right.  \mbox{ for }i,j = 1,\ldots, m$}
  \STATE{ $\dot{U} = \left(\dot{Y} V + U(\Sigma \Gamma - \dot{\Sigma}) \right) \Sigma^{-1}.$}
  \ENSURE{$\dot{U}, \dot\Sigma=\diag(\dot\sigma_1,\ldots,\dot\sigma_m), \dot{V}$}
\end{algorithmic}
\end{algorithm}
The above algorithm is mathematical `folklore', a proof can be found in, e.g.,  \cite{HayBorggaardPelletier2009}.
Note that $U^T\dot{U}$ with $\dot U$ as above is skew-symmetric, so that indeed $\dot{U} \in T_{U}St(n,m)$.
The above equations make use of the inverse $\Sigma^{-1}$ and are therefore unstable, if the singular values are small.
This effect can be alleviated by truncating the SVD to the $r\leq m$ dominant singular values.
The derivative matrices for the truncated SVD are stated in Alg. \ref{alg:SVDtruncdiff}.
\begin{algorithm}
\caption{Differentiating the truncated SVD}
\label{alg:SVDtruncdiff}
\begin{algorithmic}[1]
  \REQUIRE{matrices $Y, \dot{Y}\in\R^{n\times m}$, (truncated) SVD $Y \approx U_r\Sigma_r V^T_r$
  with $U_r\in St(n,r)$, $\Sigma_r\in \R^{r\times r}$, $V = (V_r,V_{m-r})\in O(m)$, $r\leq m\leq n$.}
  \STATE{ $\dot{\sigma}_j = (u_j)^T \dot{Y} v_j\mbox{ for }j = 1,\ldots, r$}
  \STATE{ $\dot{V}_r =V \Gamma,\quad  \Gamma_{ij} =
            \left\{
                \begin{array}{lll}
                 \frac{ \sigma_i (u_i^T \dot{Y} v_j) +  \sigma_j(u_j^T\dot{Y} v_i) }{(\sigma_j + \sigma_i)(\sigma_j - \sigma_i)}, 
                 & i\neq j,& i = 1,\ldots, m, j = 1,\ldots, r\\
                0,          & i=j, &i,j= 1,\ldots, r
                \end{array}
            \right.$ 
        } 
        \COMMENT{{//}$\quad \Gamma = \begin{pmatrix}
                             \Gamma_r\\
                             \Gamma_{m-r}
                            \end{pmatrix}
        \in \R^{m\times r}$}
  \STATE{ $\dot{U}_r = \left(\dot{Y} V_r + U_r(\Sigma_r \Gamma_r - \dot{\Sigma}_r) \right) \Sigma^{-1}_r.$}
  \ENSURE{$\dot{U}_r\in T_USt(n,r), \dot\Sigma_r=\diag(\dot\sigma_1,\ldots,\dot\sigma_r), \dot{V}_r\in T_USt(m,r)$}
\end{algorithmic}
\end{algorithm}
Since this algorithm is based on representing the derivative vectors $\dot v_j$ in terms of an eigenvector ONB
$V = (V_r, V_{m-r}) = (v_1,\ldots,v_r,v_{r+1},v_m)$, a full square orthogonal $V$ is required also in the truncated case.
Yet, note that the columns of $V_{m-r}$ feature only in the equation for $\dot V_r = V\Gamma = V_r\Gamma_r + V_{m-r}\Gamma_{m-r}$
while all other quantities are independent of $V_{m-r}$.

If the rank of $Y\in\R^{n\times m}$ is exactly $r\leq m$ and is fixed for all $t$, then the computation of the entries of the lower block $\Gamma_{m-r}$
reduces to
$\Gamma_{ij} = \frac{ u_j^T\dot{Y} v_i }{\sigma_j},  i = r+1,\ldots, m, j = 1,\ldots,r$.
In this case, the singular value matrix features a lower-right zero diagonal block
$\Sigma = \diag(\sigma_1,\ldots,\sigma_r,\sigma_{r+1},\ldots, \sigma_m)$.
In general, computing the derivatives in the presence of multiple singular values/eigenvalues is sophisticated
\cite{Alekseevsky1998}. Here, however, it is sufficient to compute the singular vectors 
$V_r(t) = (v_1(t),\ldots,v_r(t))$ associated with the pairwise distinct singular values 
and to perform a $t$-dependent orthogonal completion $V(t) =  (v_1(t),\ldots,v_r(t), v_{r+1}(t),\ldots,v_m(t))$
via the modified Gram-Schmidt process, which is differentiable.

\section{The Riemannian Stiefel log algorithm}
\label{supp:algStLog}
All numerical experiments featured in this work are performed with a SciPy\cite{SciPy2001} implementation of the following algorithm,
for the details, see \cite{StiefelLog_Zimmermann2017}.
\begin{algorithm}
\caption{Stiefel logarithm}
\label{alg:Stlog}
\begin{algorithmic}[1]
  \REQUIRE{base point $U\in St(n,p)$ and $\tilde{U}\in St(n,p)$ `close' to base point, $\tau>0$
  convergence threshold}
  \STATE{ $M:= U^T\tilde{U} \in \R^{p\times p}$}
  \STATE{ $QN := \tilde{U} - UM \in \R^{n \times p}$}
    \hfill \COMMENT{(thin) qr-decomp. of normal component of $\tilde{U}$}
  \STATE{ $V_0 := \begin{pmatrix}M&X_0\\N&Y_0\end{pmatrix} \in O_{2p\times 2p}$}
    \hfill \COMMENT{orthogonal completion}

  \FOR{ $k=0,1,2,\ldots$}
  \STATE{$\begin{pmatrix}A_k & -B_k^T\\ B_k & C_k\end{pmatrix} := \log_m(V_k)$}
    \hfill \COMMENT{matrix log, $A_k, C_k$ skew}
  \IF{$\|C_k\|_F \leq \tau$}
    \STATE{break}
  \ENDIF
  \STATE{ $\Phi_{k} = \exp_m{(-C_k)}$}
    \hfill \COMMENT{matrix exp, $\Phi_{k}$ orthogonal}
  \STATE{$V_{k+1} := V_{k}W_k$, where $W_k:=\begin{pmatrix}I_p& 0\\ 0 & \Phi_k\end{pmatrix}$
  }
    \hfill \COMMENT{update}
  \ENDFOR
  \ENSURE{$\Delta := Log_{U}^{St}(\tilde{U}) = U A_k + QB_k \in T_{U}St(n,p)$}
%
\end{algorithmic}
\end{algorithm}
%

\end{document}